\documentclass[12pt]{article}
\usepackage{amssymb}
\usepackage{verbatim}

\usepackage[latin1]{inputenc}
\def\bea{\begin{eqnarray*} }
\def\eea{\end{eqnarray*} }
\newtheorem{definition}{Definition}
\newtheorem{example}[definition]{Example}
\newtheorem{proposition}[definition]{Proposition}
\newtheorem{lemma}[definition]{Lemma}
\newtheorem{corollary}[definition]{Corollary}
\newtheorem{theorem}[definition]{Theorem}
\newtheorem{remark}[definition]{Remark}
\newcommand{\hess}{\mathrm{Hess}}
\newcommand{\niu}{\mbox{{\Large $\nu$}}}
\newcommand{\rank}{\mathrm{rank}}
\newcommand{\po}{{\hspace*{-1ex}}{\bf .  }}
\newcommand{\R}{\mathbb{R}}
\newcommand{\Q}{\mathbb{Q}}
\newcommand{\Sf}{\mathbb{S}}
\newcommand{\spa}{\mbox{span}}

\def\e{\epsilon}

\newcommand{\be}{\begin{equation} }
\newcommand{\ee}{\end{equation} }
\newcommand{\Ker}{\mbox{Ker} }
\def\<{\langle}
\def\>{\rangle}
\def\proof{\noindent{\it Proof: }}
\def\qed{\ifhmode\unskip\nobreak\fi\ifmmode\ifinner\else
\hskip5 pt \fi\fi\hbox{\hskip5 pt \vrule width4 pt
height6 pt  depth1.5 pt \hskip 1pt }}

\begin{document}

\title{Parabolic submanifolds of rank two}
\author {Marcos Dajczer \& Pedro Morais}
\date{}
\maketitle

An immersed submanifold $f\colon M^n\rightarrow\R^N$, $n\ge 3$,  into Euclidean space with the induced metric is called of \emph{rank two} if at any point the
kernel of its vector valued second fundamental form has codimension two. Equivalently, we have that the image of the Gauss map in the Grassmannian of 
non-oriented $n$-planes $G_n^N$ is a surface. These submanifolds have been the object of a great deal of work in Riemannian Geometry since long time ago.  For instance, see \cite {BKV} and references therein. This interest is in good part motivated  
by the fact that their curvature tensor is ``as flat as possible" without vanishing altogether.

The subspace  spanned by the second fundamental form, usually called the first normal space and denoted by $N_1$, of a rank two submanifold satisfies $\dim N_1\leq 3$  at any point. It turns out that if in substantial codimension, any  rank two submanifold is a hypersurface if $\dim N_1=1$ at any point. Then $f$ is either a Euclidean surface or the cone over a spherical surface, up to a Euclidean factor, if dim $N_1=3$ everywhere. Submanifolds in the remaining and much more interesting case, namely,  when $\dim N_1=2$ everywhere, have been divided in three classes: elliptic, hyperbolic and parabolic. A complete parametric description of the elliptic submanifolds  was given in \cite{DajFlo1}.

For codimension $N-n=2$, it was shown in \cite{DajFlo2} that elliptic and nonruled parabolic submanifolds are \emph{genuinely} rigid. This means that given any other isometric immersion $\tilde f\colon\,M^n\to\R^{n+2}$ there is an open dense  subset of $M^n$ such that
restricted to any connected component $f|_U$ and $\tilde{f}|_U$ are either congruent or there are an isometric
embedding \mbox{$j\colon\, U\hookrightarrow N^{n+1}$} into a Riemannian manifold $N^{n+1}$ and either flat or
isometric noncongruent hypersurfaces $F,\tilde F\colon\,N^{n+1}\to\R^{n+2}$\ such that  $f|_U=F\circ j$ and
$\tilde{f}|_U=\tilde F\circ j$. Recently, we  proved \cite{DajMor} that nonruled parabolic submanifolds in codimension two are not only genuinely rigid but, in fact, isometrically rigid. 

The goal of this paper is to classify parametrically  parabolic submanifolds  in any codimension. First, we describe  the ones that are ruled and  show that they are the only parabolic submanifolds that admit an isometric immersion as a hypersurface. Then, we  classify  the nonruled ones by two different means.  In fact, we provide the \emph{polar} and \emph{bipolar} parametrizations, each of which is associated 
to a parabolic surface and a function on the surface which satisfies a parabolic differential equation. To conclude, we describe the structure of the singular set of the nonruled parabolic submanifolds.

\section{Parabolic submanifolds.}

In this section, we introduce the concept of parabolic submanifold and  study in detail the structure of the  normal bundle.\vspace{1,5ex}

  We denote by $f\colon\,M^n\to\Q_\e^N$, $\e=0,1$, a connected  $n$-dimensional submanifold of either Euclidean space $\R^N$ ($\e=0$) or unit Euclidean sphere $\Sf^N$ ($\e =1$) with codimension $N-n$. The \mbox{{\it $k^{th}$-normal space\/}} $N^f_k(x)$ of $f$ at $x\in M^n$ is defined as
$$
N^f_k(x)=\spa\{\alpha_f^{k+1}(X_1,\ldots,X_{k+1})\,;X_1,\ldots,X_{k+1}\in T_xM\}.
$$
Here, $\alpha_f^\ell\,\colon\,TM\times\cdots\times TM\to T_f^\perp M$, $\ell\geq 2$, is the symmetric tensor known as 
the~\emph{$\ell^{th}$-fundamental form} and given by
$$
\alpha_f^\ell(X_1,\ldots,X_\ell)=\pi^{\ell-1}\left(\nabla^\perp_{X_\ell}\ldots
\nabla^\perp_{X_3}\alpha_f(X_2,X_1)\right)
$$
where $\pi^\ell$ stands for the orthogonal projection
$\pi^\ell\colon\,T_f^\perp M\to (N^f_1\oplus\ldots\oplus N^f_{\ell-1})^\perp$ 
and $T_f^\perp M$ is endowed with the normal connection $\nabla^\perp$ induced by the metric connection $\tilde\nabla$ in the ambient space. We agree that
$\alpha_f^1\colon\,TM\to TM$ is $\alpha_f^1=I$ and
denote $\alpha_f^2=\alpha_f$ ($\pi^1=I$) as usual.

We always assume that $f\colon\,M^n\to\Q_\e^N$ is substantial and has rank~$2$. The later condition is denoted as
$\rank_f=2$, and means that the relative nullity subspaces $\Delta(x)\subset T_xM$ defined as
$$
\Delta(x)=\{X\in T_xM : \alpha_f(X,Y)=0 \,;\, Y\in T_xM\},
$$
form a tangent subbundle of codimension two. It is a standard fact  that the relative nullity distribution is integrable and that the leaves are totally geodesic submanifolds of the ambient space $\Q_\e^N$.

  The cone \mbox{$Cf\colon\;M^n\times\R_+\to\R^{N+1}$} of a submanifold   \mbox{$f\colon\;M^n\to\Sf^N$} of rank two has the same rank
since the relative nullity leaves of $Cf$ are the cones of the relative nullity leaves of $f$. Moreover, one has
that $N^{Cf}_k=~N^f_k$, $k\geq 1$, up to parallel transport in $\R^{N+1}$. Thus, it suffices to consider the
Euclidean case since we had restricted ourselves to submanifolds of $\R^N$ \mbox{and $\Sf^N$}.

 The condition $\rank_f=2$ and the symmetry of the second fundamental form  imply that the first normal spaces of $f$ satisfy $\dim N^f_1\leq 3$ at any point. By Theorem~$1$ in \cite{DajToj} we have that $f$ is a hypersurface in substantial codimension if
\mbox{$\dim N^f_1=1$} everywhere. On the other hand, it is not difficult to show that a  submanifold with $\dim N^f_1=3$
everywhere is either a Euclidean surface  or the cone over a spherical surface up to Euclidean factor.  In the
remaining case  when $\dim N^f_1=2$ everywhere,  either there exists a pair of linearly independent
``conjugate directions" $X_1,X_2\in \Delta^\perp$ , i.e., \be\label{sign} 
\alpha_f(X_1,X_1)\pm \alpha_f(X_2,X_2)=0, 
\ee 
or $f$ admits an ``asymptotic direction" $0\neq Z\in\Delta^\perp$, i.e.,
$
\alpha_f(Z,Z)=0.
$
In cases (\ref{sign}) the submanifold was called  \emph{elliptic} for the plus sign and \emph{hyperbolic} for the minus sign in \cite{DajFlo1}.
\begin{definition}\po {\em A submanifold  $f\colon\,M^n\rightarrow \Q^N_\epsilon$ is called  \emph{parabolic} if we have:
\begin{itemize}
\item [(i)]$\rank_f =2$,
\item [(ii)]$\dim N^f_1=2$,
\item [(iii)] There is a nonsingular asymptotic vector field
$Z\in\Delta^\perp$, i.e., $\alpha_f (Z,Z)=0$.
\end{itemize}
}\end{definition}

Notice that cones of parabolic spherical submanifolds are also parabolic.
\medskip

Let $f\colon\,M^n \rightarrow \R^N$ be a parabolic submanifold. We always denote by $\{X,Z\}$ an orthonormal
frame in $\Delta^{\perp}$ where  $Z$ is an asymptotic vector field. Clearly, we can always take an orthonormal smooth frame
$\{\eta_1,\eta_2\}$ in $N_1^f$ such that the shape operators take \
the form \be\label{forma}
A^f_{\eta_1}|_{\Delta^\perp}=\left[\begin{array}{cc}
  a & b \\
  b & 0
\end{array}\right]\quad\textrm{and}\quad A^f_{\eta_2}|_{\Delta^\perp}=\left[\begin{array}{cc}
  c & 0 \\
  0 & 0
\end{array}\right]
\ee
where the functions $b,c$ never vanish. In particular, we see that 
the asymptotic field $Z$ is unique up to sign.\vspace{1ex}

An easy argument given in \cite{DajFlo1} proves the following fact.

\begin{proposition}\po\label{norm}  Assume that $f\colon\,M^n\to\Q_\e^N$  satisfies
$\dim N_1^f=2$ at any point. Then, we have that  
$\dim N^f_{k}\leq~2$ for all $k\geq 1$.
\end{proposition}

We always admit that the fibers of any  $N^f_k$ have constant 
dimension and thus form subbundles of the normal
bundle. If $\tau = \tau^f$ denotes the index of the ``last" 
of the normal subbundles of $f$, then
$T_f^\perp M=N^f_1\oplus\cdots\oplus N^f_\tau$ 
since, by assumption, $f$ is substantial. \vspace{1ex}

We denote
$$
\xi_1^k = \alpha_f^{k+1}(X , \ldots ,X)
\quad \textrm{and}\quad  \xi_2^k
= \alpha_f^{k+1}(Z,X, \ldots , X).
$$
Since $\alpha_f^{k+1}(Z,Z,Y_1,\ldots,Y_{k-1})=0$, it is clear that
$$
N_k^f=\spa\{\xi_1^k, \xi_2^k\}\;\;\;\mbox{for} \;\;\;1\le k\le\tau^f. 
$$

\begin{proposition}\po\label{Lemma1}
For $\; 1\leq k\leq \tau^f-1$ the following holds:
\begin{itemize}
\item[(i)] $(\tilde{\nabla}_Z \: \xi_1^k)_{N^f_{k+1}}=(\tilde{\nabla}_X
\:\xi_2^k)_{N^f_{k+1}}=\xi_2^{k+1},$
\item[(ii)] $(\tilde{\nabla}_X \: \xi_1^k)_{N^f_{k+1}}=\xi_1^{k+1},$
\item[(iii)]$(\tilde{\nabla}_Z \:\xi_2^k)_{N^f_{k+1}}=0$.
\end{itemize}
\end{proposition}

\proof From the definition of the  k-normal spaces, given $\eta\in N^f_l$ we have 
\be\label{comp}
\nabla^{\perp}_Y\eta \in N^f_{l-1}\oplus N^f_l \oplus N^f_{l+1} 
\ee 
where  $N^f_0=0=N^f_{\tau^f+1}$. Then, 
$$
\xi_2^{k+1}
=(\nabla^{\perp}_Z(\nabla^{\perp}_X
\ldots\nabla^{\perp}_X\alpha_f(X,X))_{N^f_k})_{N^f_{k+1}}
=(\tilde{\nabla}_Z \: \xi_1^k)_{N^f_{k+1}},
$$
$$
\xi_2^{k+1}
=(\nabla^{\perp}_X(\nabla^{\perp}_Z
\ldots\nabla^{\perp}_X\alpha_f(X,X))_{N^f_k})_{N^f_{k+1}}
=(\tilde{\nabla}_X \: \xi_2^k)_{N^f_{k+1}}
$$
and $(i)$ has been proved. The proof of $(ii)$ is similar. 
For $(iii)$, we have
$$
(\tilde{\nabla}_Z\xi_2^k)_{N^f_{k+1}}
=(\nabla^{\perp}_Z(\nabla^{\perp}_X\ldots
\nabla^{\perp}_X\alpha_f(X,Z))_{N^f_k})_{N^f_{k+1}}
=\alpha_f^{k+2}(X,\ldots,Z,Z)
=0.\qed
$$

The following fact was proved in \cite{DajFlo1}.

\begin{proposition}\po\label{nparalelos} If $f\colon\,M^n\rightarrow\R^N$ is a parabolic submanifold, then the normal subbundles $N_k^f,\; 1\leq k\leq \tau^f$, are parallel
in $\R^N$ along $\Delta$.
\end{proposition}

   Let  $\niu_k\subset N^f_k\times N^f_k,\; 0\leq k\leq\tau^f$, 
be the subspace defined as
$$
\displaystyle \niu_k=\{ (\mu_1, \mu_2) \in N^f_k \times
N^f_k : \< \mu_2, \xi_2^k \> =0 \;\; \textrm{and}\;\; \< \mu_2, \xi_1^k\> = \< \mu_1,
\xi_2^k\> \}.
$$
It is easy to see that $\niu_k$ is independent  of the base 
$\{X,Z\}$ with $Z$ asymptotic. Clearly,  $\xi_1^k=0$ implies 
that $\niu_k=0$. We also have the following facts.

\begin{lemma}\po\label{Lemma:proj} For $\; 1\leq k\leq \tau^f$ 
the following holds:
\begin{itemize}
\item[(i)]  $\displaystyle \dim \niu_k=2 
\;\mbox{if and only if}\; \dim N^f_k=2$,
\item[(ii)]  $\displaystyle \dim\niu_k=1\; 
\mbox{if and only if} \;\dim N^f_k=1 \;\;
\mbox{and} \;\;\xi_2^k = 0$,
\item[(iii)]  $\displaystyle\dim\niu_k=0 \; 
\textrm{if and only if} \;\dim N^f_k=1 \;\;\mbox{and} \;\;
\xi_2^k \neq 0$.
\end{itemize}
\end{lemma}

\proof If $\dim \niu_k=2$, we either may choose
$(\mu_1,\mu_2)\in \niu_k$ such that $\mu_1\neq 0\neq\mu_2$ 
or we are done. It is easy to see that $\mu_1$ and $\mu_2$ 
must be linearly independent, and thus $\dim N_k^f=2$. Then,
take $0\neq v\in N^f_k$ such that
$\< v, \xi_2^k \> =0$, and set  
$u=(\< v, \xi_1^k \>/\|\xi^k_2\|^2)\xi^k_2$.
Hence,  $u,v$ are a base of $N^f_k$ and 
$(u,v),(u+v,v) \in \niu_k $ are linearly independent. 
This proves  $(i)$.
The proofs of $(ii)$ and $(iii)$ follow easily form the  definition~of $\niu_k$. \qed

\begin{definition}\po
{\em Given a parabolic submanifold $f\colon\,M^n\to\Q_\e^N\subseteq\R^{N+\e}$, we call an element $\beta\in C^\infty(M^n,\R^{N+\e})$ a {\it $k$--cross section\/}  
to $f$, $1\leq k\leq\tau^f$, if at any point
$$
\beta_\ast\,(TM)\subset N_{k}^f\oplus\cdots\oplus N^f_{\tau^f},
$$
up to parallel transport in $\R^{N+\e}$.}
\end{definition}

\begin{lemma}\po\label{Lemma:proj1} Let
$\mathcal{P}_k\colon\,C^{\infty}(M^n,\R^{N+\epsilon})\rightarrow N^f_k \times N^f_k$, $1\leq k \leq \tau^f$, be the tensor
$$
\mathcal{P}_k(\beta)=((\beta_{\ast}X)_{N^f_k}, (\beta_{\ast}Z)_{N^f_k}).
$$
Then $\mathcal{P}_k(\beta)\in\niu_k$  for any  $k$--cross section $\beta$
to $f$. Moreover, the tensor
$$
\mathcal{P}_k \vert_{N^f_{k+1}}\colon\,{ N^f_{k+1}} \rightarrow \niu_k, \;\;\;\; 1\leq k\leq\tau^f-1,
$$
is injective.
\end{lemma}

\proof We have,
\begin{eqnarray*}
\<\beta_*X,\xi_2^k \>
\!\!&=&\!\! \< \tilde{\nabla}_{X}
\beta,\tilde{\nabla}_{Z}(\nabla^\perp_{X}
\ldots\nabla^\perp_{X}\alpha_f(X,X))\>\\
\!\!&=&\!\! Z\< \tilde{\nabla}_{X} \beta,
\nabla^\perp_{X}\ldots\nabla^\perp_{X}\alpha_f(X,X)\>
-\< \tilde{\nabla}_{Z}\tilde{\nabla}_{X} \beta,\nabla^\perp_{X}
\ldots\nabla^\perp_{X}\alpha_f(X,X)\>\\
\!\!&=&\!\! \<
\tilde{\nabla}_{Z}\beta,
\tilde{\nabla}_{X}(\nabla^\perp_{X}\ldots\nabla^\perp_{X}\alpha_f(X,X))\>\\
\!\!&=&\!\! \<\beta_*Z,\xi_1^k \>.
\end{eqnarray*}
A similarly argument gives
$$
\<\beta_*Z, \xi_2^k \>
= \<\beta_*X,\alpha_f^{k+1}(Z,Z,X,\ldots,X)\>=0.
$$
To conclude, observe that if
$\eta\in N^f_{k+1}$ satisfies $\mathcal{P}_k (\eta)=0 $, then
$$
0=\<\eta_*X,\xi_j^k\>=\<\tilde{\nabla}_X\eta,\xi_j^k\>=
-\<\eta, \tilde{\nabla}_X\xi_j^{k}\>=
-\< \eta,\xi_j^{k+1} \>, \;\: j=1,2.
$$
Hence, $\eta=0$.
\qed

\begin{proposition}\po\label{proj2}
Let $f\colon\,M^n\rightarrow\R^N$ be a parabolic submanifold. Then, we have:
\begin{itemize}
\item[(i)] $\xi_1^k\neq 0$  for any $1\leq k\leq\tau^f-1$,
\item[(ii)] $\xi_2^k=0$ if and only if $\dim N_k^f=1$,
\item[(iii)] If $\xi_2^k=0$, then $\xi_2^j=0$ for $j\geq k$.
\end{itemize}

\end{proposition}

\proof To prove $(i)$ suppose that $\xi_1^k=0$. Thus, $\niu_k=0$. Then Lemma \ref{Lemma:proj1} gives $N^f_{k+1}=0$,
which is not possible. For  $(ii)$ suppose that $\dim N_k^f=1$ and $\xi_2^k\neq 0$. We have that
$\niu_k=0$ from Lemma \ref{Lemma:proj}, and by Lemma \ref{Lemma:proj1} this is a contradiction. Finally, to prove $(iii)$ assume $\xi_2^k=0$.
Using (\ref{comp}) we have
\begin{eqnarray*}
\xi_2^{k+1}
\!\!&=&\!\!\pi^{k+1}( \nabla^\perp_X\nabla^\perp_Z\nabla^\perp_X
\ldots\nabla^\perp_X\alpha_f(X,X))\\
\!\!&=&\!\!\pi^{k+1} (\nabla^\perp_X(\pi^k (\nabla^\perp_Z\nabla^\perp_X
\ldots\nabla^\perp_X\alpha_f(X,X)))\\
\!\!&=&\!\!\pi^{k+1} (\nabla^\perp_X\xi_2^k)=0.\qed
\end{eqnarray*}

\begin{definition}\po {\em We say that a parabolic submanifold
$f\colon\,M^n \rightarrow \Q^N_{\epsilon}$ has
\emph{critical index}
$\tau_0^f \in \{1,\ldots,\tau^f-1\}$ if $\xi_2^{\tau_0^f}\neq 0$ and $\xi_2^k=0$ for any $k\geq\tau_0^f+1$.}
\end{definition}

\begin{corollary}\po\label{Proposition:isom}
Assume that $f$ possesses critical index. Then:
\begin{itemize}
\item[(i)] $ \dim N^f_k=2,\; 1\leq k \leq \tau_0^f,$
\item[(ii)] $ \dim N^f_k=1,\; \tau_0^f+1 \leq k \leq \tau^f,$
\item[(iii)] The tensor, $\mathcal{P}_k \vert_{ N^f_{k+1}}\colon\,N^f_{k+1} \rightarrow \niu_k$ is an isomorphism for $k\leq\tau_0^f-1$.
\end{itemize}
\end{corollary}

\section{Intrinsic proprieties}

In this section we analyze the  metric structure of the parabolic submanifolds.

\begin{proposition}\po\label{prop:parabolica}  Let $f\colon\,M^n \rightarrow \R^N$ be a  parabolic submanifold. Then, 
$$
\mathcal{F}=\spa\{Z\}\oplus\Delta
$$
is an integrable distribution and the leaves are flat hypersurfaces.
\end{proposition}

\proof  We first show that the line bundle  $L=\spa\{\xi_2^1\}$ is parallel along the leaves of relative nullity.
The unit vector field $\eta\in N_1^f $ orthogonal to $\xi_2^1$ is the only one, up to sign, such that $A^f_\eta
Z=0$. Thus  $A^f_\eta$ has rank $1$. In view of Proposition \ref{nparalelos} it is sufficient to show that $\eta$
is parallel along $\Delta$.

Recall that the \emph{splitting tensor} $C$ associates to
$T\in\Delta$ the endomorphism $C_T$ of $\Delta^\perp$ defined  as
$$
C_T X=-\left(\nabla_X T\right)_{\Delta^\perp}.
$$
It is well-known \cite{DajGro} that  the differential equation
\be\label{cete}
\nabla _T A^f_\xi=A^f_\xi \circ C_T
\ee
is satisfied along $\Delta^\perp$ if $\xi\in T_f^\perp M$ is parallel along $\Delta$.

Let $x\in M^n$ and $\gamma$ a geodesic with $\gamma(0)=x$ contained  in the corresponding leaf of $\Delta$. If
$\delta_t$ is the parallel transport of $\eta_x$ along $\gamma$, we have
$$
\nabla_{\gamma'}A^f_{\delta_t}=A^f_{\delta_t}\circ C_{\gamma'}.
$$
Hence,
$A^f_{\delta_t}=A^f_{\eta_x}e^{\,\int_0^t C_{\gamma'}d\tau}$.
Thus  $A^f_{\delta_t}$ has rank $1$ and, therefore $\eta=\delta_t$ is parallel.

Since the left hand side of
$$
\nabla_T A_\eta^f= A_\eta^f\circ C_T
$$
is symmetric, we obtain that
$$
A_\eta^fC_T Z=C_T^tA_\eta^fZ=0.
$$
Thus $C_TZ\in\spa\{Z\}$, that is, $\< \nabla_Z T,X\>=0$.  Then 
the Codazzi equation yields
$$
\nabla_T^\perp\alpha_f(Z,X)
-\<\nabla_TZ,X\>\alpha_f(X,X)
+\<\nabla_ZT,Z\>\alpha_f(Z,X)=0.
$$
Using that $L$ is parallel along  $\Delta$,  we obtain that
$\<\nabla_TZ,X\>=0$. Hence 
$\mathcal{F}$ is integrable. Moreover, the second fundamental form 
of a leaf $U$ is
$$
A^U_X=\left[\begin{array}{cc}
  \lambda & 0 \\
  0 & 0
\end{array}\right]
$$
where $\lambda=\<\nabla_Z Z,X\>$. Thus the leaves of $\mathcal{F}$ are flat. \qed\vspace{1,5ex}

Recall that a submanifold $f\colon\,M^n\rightarrow \Q^N_\epsilon$ is called \emph{ruled}  when $M^n$ admits a hypersurface foliation of totally geodesic submanifolds of $\Q^N_\epsilon$.

\begin{example}\po\label{obs:sup} {\em Ruled Euclidean submanifolds of rank $2$ without flat points and substantial codimension at least $2$ are basic examples of parabolic submanifolds. In fact, it follows from Corollary 4.7 in \cite{Daj} that $\dim N_1^f=2$.
}\end{example}

  From the proof of Proposition \ref{prop:parabolica}  we have the following fact.

\begin{corollary}\po Let $f\colon\,M^n\rightarrow \R^{N}$ be a ruled parabolic submanifold. Then the leaves of $\mathcal{F}$ are totally geodesic in $M^n$.
\end{corollary}

\section{Regularity}

A key ingredient in the parametric description of the elliptic submanifolds given in  \cite{DajFlo1} was the
regularity of the  $k$-normal spaces. In fact,  any elliptic subma\-nifold  $f$  satisfies $\dim N_k^f=2,$ $ 1\leq
k\leq \tau^f-1 $, whereas the dimension of $N_{\tau^f}^f$ is determined by the codimension. In this paper, that a
parabolic submanifold is regular roughly means that the $N_k^f$'s behave as in the elliptic case.  The main result in
this section is that nonregular parabolic submanifolds are necessarily ruled.

\begin{definition}\po{\em We say that a parabolic submanifold  $f\colon\,M^n\!\rightarrow \R^N$ is \emph{regular} if
$\dim N_k^f=2$ for any $1\leq k\leq \tau^f-1$.}
\end{definition}

By Corollary \ref{Proposition:isom}, the following holds:
$$
f \mbox{ is regular if and only if } \left\{\begin{array}{l} \!\dim N^f_{\tau^f}=2 \;\;\iff
\;\xi_2^{\tau^f}\neq 0,\;\mbox{if}\;\; N-n \quad \mbox{is even}
\vspace*{1ex}\\
\!\dim N^f_{\tau^f-1}=2 \!\!\iff \xi_2^{\tau^f-1}\neq 0,\;\!\mbox{if}\;\;
 N-n \quad \mbox{is odd}.
\end{array}
\right.
$$ 
Observe that ruled surfaces with $\dim N_1=2$ are parabolic. We give next an example of such a surface that is nonregular.

\begin{example}\po{\em Let
$c\colon\,I\subset\R \rightarrow \R^6 $ be a smooth curve parametrized by arc length with Frenet frame $E_1,\ldots, E_6$ and constant Frenet curvatures $k_j\neq 0, 1\leq j \leq 5$.
The map $X\colon\,\R^2 \rightarrow \R^6$ given by
$$
X(s,t)=c(s)+tE_2(s)
$$
parametrizes a substantial complete surface that is parabolic for $t\neq 0$.
An easy calculation gives $\xi_2^2=0$, that is, $\tau_0^X=1$. Hence, $\dim N_2^{X}=1$ and therefore $X$ is nonregular.
}\end{example}

By a parabolic submanifold being {\it nonruled\/} we understand that none of the leaves of $\mathcal{F}$ is
totally geodesic in $M^n$ or, equivalently, in $\R^N$.

\begin{theorem}\po\label{main}
Nonruled parabolic submanifolds  $f\colon\,M^n\rightarrow \R^N$ are  regular.
\end{theorem}

The proof of Theorem \ref{main} will follow from two results. First, we give a sufficient condition for a parabolic submanifold in odd codimension to be ruled.

\begin{proposition}\po\label{regrada} Let $f\colon\,M^n\rightarrow \R^N$ be a regular parabolic submanifold satisfying that $\xi_2^{\tau^f}=0$ at any point.
Then $f$ is ruled.
\end{proposition}

\proof We claim that $f$ is ruled if and only if  $L=\spa\{\xi_2^1\}$ is parallel along $\mathcal{F}$. From the
proof of Proposition \ref{prop:parabolica}, we know that $L$ is parallel along $\Delta$. Clearly, that $f$ is ruled
is equivalent to $\nabla_ZZ=0$. Take an orthonormal frame $\{\eta_1,\eta_2\}$ in $N_1^f$ as in (\ref{forma}).
Since $\eta_1\in L$, we have to show that 
\be\label{normalregrada} 
\nabla_ZZ=0\;\; \mbox{if and only
if}\;\;(\nabla^\perp_Z\eta_1)_{N_1^f}=0. 
\ee 
From the Codazzi equation 
$$ 
\<(\nabla_XA^f_{\eta_2})Z
-(\nabla_ZA^f_{\eta_2})X,Z\>=0, 
$$ 
we get
$$
c\< \nabla_ZZ,X\>
=b\<\nabla^\perp_Z\eta_1,\eta_2\>.
$$
Being $f$ parabolic we obtain  $b\neq 0\neq c$, and the claim follows. \vspace{1ex}

We first consider the case $N-n=3$. We have, $\dim N_1^f=2,\; \dim N_2^f=1$ and $\xi_2^2=0$. It suffices to show that $\eta_1$ is parallel along  $Z$. By Proposition \ref{Lemma1}, the subbundles $N_1^f, N_2^f$ are parallel along  $Z$. Thus, the Codazzi equation gives
$$
A^f_{\nabla_X^\perp\delta}Z=A^f_{\nabla_Z^\perp\delta}X=0
$$
where $\delta\in N_2^f$ has unit length. Using (\ref{forma}) we obtain 
\be\label{projn1}
\left(\nabla^\perp_X\delta\right)_{N_1^f}\in \spa\{\eta_2\}. 
\ee 
From $X\< \eta_1, \delta\>=0$  and
(\ref{projn1}) we have 
\be\label{eta1} 
\left(\nabla^\perp_X\eta_1\right)_{N^f_2}=0. 
\ee 
The Ricci equation, using (\ref{projn1}), (\ref{eta1}) and the  parallelism of $N_1^f$ along  $Z$ gives
\begin{eqnarray*}
0=\< R^\perp(X,Z)\eta_1,\delta\>\!\!&=&\!\!\< \nabla^\perp_X \nabla^\perp _Z \eta_1-\nabla^\perp_Z \nabla^\perp _X
\eta_1-\nabla^\perp_{[X,Z]}\eta_1,  \delta\>\\
\!\!&=&\!\!\<\nabla^\perp_X\nabla^\perp _Z \eta_1,
\delta\>=-\<\nabla^\perp_Z\eta_1,\nabla^\perp_X\delta\>\\
 \!\!&=&\!\!\<\nabla^\perp_Z\eta_1,\eta_2\>\<\nabla^\perp_X \eta_2,\delta\>.
\end{eqnarray*}
But $\<\nabla^\perp_X \eta_2,\delta\>\neq 0$ since $N_1^f$ is not parallel. Thus,
$(\nabla^\perp_Z\eta_1)_{N_1^f}=0$.
\vspace{1ex}

  We now consider the general case $N-n\ge 5$. Take an orthonormal basis $\{\eta^k_1,\eta^k_2\}$ of $N_k^f$ for any $1\leq k\leq \tau^f-1$ such that
$$
\xi_1^k=a_k\eta^k_1+c_k\eta^k_2\quad\textrm{and}\quad \xi_2^k=b_k\eta^k_1.
$$
Proposition \ref{Lemma1} gives 
\be\label{igualcompG} \left(\nabla^\perp_Z\eta^{k}_1\right)_{N_{k+1}^f}=0
\quad\textrm{and}\quad c_{k}\left(\nabla^\perp_Z\eta^{k}_2\right)_{N_{k+1}^f}=
b_{k}\left(\nabla^\perp_X\eta^{k}_1\right)_{N_{k+1}^f}. 
\ee 
Since  $\dim N_{k}^f=2$, $1\leq k\leq \tau^f-1$, it
follows from (\ref{igualcompG}) that 
\be\label{baseG}
N_{k}^f=\spa\left\{\left(\nabla^\perp_X\eta^{k-1}_1\right)_{N_{k}^f},
\left(\nabla^\perp_X\eta^{k-1}_2\right)_{N_{k}^f}\right\}. 
\ee 
From  (\ref{igualcompG}) and $\xi_2^{\tau^f}=0$, we have
\be\label{compn3G} 
(\nabla^\perp_Z\eta_1^{\tau^f-1})_{N_{\tau^f}^f}
=(\nabla^\perp_X\eta^{\tau^f-1}_1)_{N_{\tau^f}^f}= (\nabla^\perp_Z\eta^{\tau^f-1}_2)_{N^f_{\tau^f}}=0. 
\ee 
Thus 
$N_1^f \oplus\ldots\oplus N_{\tau^f-1}^f$ and $N_{\tau^f}^f$ 
are both parallel along $Z$.  The Ricci equation for
$\delta \in N_{\tau^f}^f$  and (\ref{compn3G}) give
\begin{eqnarray*}
0=\< R^\perp(X,Z)\eta_1^{\tau^f-1},\delta\>
\!\!&=&\!\!\< \nabla^\perp_X \nabla^\perp_Z \eta_1^{\tau^f-1}, \delta\>
=-\<\nabla^\perp_Z\eta_1^{\tau^f-1},
\nabla^\perp_X\delta
\>\\
\!\!&=&\!\!\< \nabla^\perp_Z \eta_1^{\tau^f-1},
\eta_2^{\tau^f-1}\> \<
\nabla^\perp_X \eta_2^{\tau^f-1},\delta\>.
\end{eqnarray*}
But $\<\nabla^\perp_X\eta_2^{\tau^f-1},\delta\>\neq 0$ since $f$ is substantial. Therefore,
$$
(\nabla^\perp_Z \eta_1^{\tau^f-1})_{N_{\tau^f-1}^f}=0.
$$
To conclude again that $\<\nabla^\perp_Z\eta_1^{1},\eta_2^{1}\>=0$, it suffices to show that if
\be\label{eta22G}
(\nabla^\perp_Z \eta_1^{\ell+1})_{N_{\ell+1}^f}=0,\;\;\;  1\le\ell\le\tau^f-2,
\ee
then
\be\label{fim}
(\nabla^\perp_Z \eta_1^{\ell})_{N^f_\ell}=0.
\ee
Being $\eta_1^\ell$ collinear with $\xi_2^\ell$ and $\eta_1^{\ell+1}$
with $\xi_2^{\ell+1}$,  then  $\eta_1^{\ell+1}$ and $(\nabla_X^\perp\eta_1^\ell)_{N_{\ell+1}^f}$ are also
collinear. From (\ref{eta22G}), we have
\be\label{teo31}
\< \nabla^\perp_Z(\nabla_X^\perp
\eta_1^\ell)_{N_{\ell+1}^f}, \eta_2^{\ell+1}\>=0.
\ee
The Ricci equation
using (\ref{igualcompG}) and (\ref{teo31}) yields
\begin{eqnarray*}
0\!\!&=&\!\!\!\!\<R^\perp(X,Z)\eta_1^\ell,\eta_2^{\ell+1}\>=\< \nabla^\perp_X \nabla^\perp _Z \eta_1^\ell-\nabla^\perp_Z \nabla^\perp _X
\eta_1^\ell
-\nabla^\perp_{[X,Z]}\eta_1^\ell,\eta_2^{\ell+1}\>\\
\!\!&=&\!\!\!\< \nabla^\perp_X\< \nabla^\perp_Z\eta_1^\ell,\eta_2^\ell\>
\eta_2^\ell,\eta_2^{\ell+1}\>
-\<\nabla^\perp_Z(\nabla^\perp_X \eta_1^\ell)_{N_\ell^f},\eta_2^{\ell+1}\>
-\<\nabla_XZ,X\>\<\nabla^\perp_X\eta_1^\ell,
\eta_2^{\ell+1}\>\\
\!\!&=&\!\!\!\<
\<\nabla^\perp_Z\eta_1^\ell,
\eta_2^\ell\> \nabla^\perp_X\eta_2^\ell
-\<\nabla^\perp_X\eta_1^\ell
\eta_2^\ell\>\nabla^\perp_Z\eta_2^\ell,
- \<\nabla_XZ,X\>\nabla^\perp_X\eta_1^\ell,
\eta_2^{\ell+1}\>.
\end{eqnarray*}
Thus,
$$
\left(\<\nabla^\perp_Z\eta_1^\ell,
\eta_2^\ell\>\nabla^\perp_X\eta_2^\ell-
\<\nabla^\perp_X\eta_1^\ell,
\eta_2^\ell\>\nabla^\perp_Z\eta_2^\ell
-\<\nabla_XZ,X\>\nabla^\perp_X\eta_1^\ell\right)_{N_{\ell+1}^f}\in\spa\{\eta_1^{\ell+1}\},
$$
and we obtain (\ref{fim}) from (\ref{igualcompG}) and (\ref{baseG}). \qed\vspace{1,5ex}

To conclude that $f$ is ruled,  from (\ref{eta22G}) and (\ref{fim}) in the proof of the preceding result it is sufficient to show that there exists an index $1\leq \ell\leq \tau^f-2$ such that $(\nabla^\perp_Z \eta_1^{\ell+1})_{N_{\ell+1}^f}=0$. Thus, this gives the following fact.

\begin{corollary}\po\label{coro:regrada}
Let $f\colon\,M^n\rightarrow \R^N$ be a regular parabolic submanifold. If there is an index $1\leq s \leq
\tau^f-1$ such that $\eta_1^s=\xi_2^s/\|\xi_2^s\| \in N^f_s$ satisfies
 $(\nabla^\perp_Z\eta_1^s)_{N^f_s}\!\!=0$, then $f$~is ruled.
\end{corollary}

Our next result deals with nonregular parabolic submanifolds.

\begin{proposition}\po\label{redcod}
Let $f\colon\,M^n\rightarrow \R^N$ be a simply connected parabolic submanifold. Assume that $\dim N_{k_0-1}^f=2$ and $\dim N_{k_0}^f =1$ for some index $2\leq k_0 \leq \tau^f-1$. Then, there exists a parabolic regular isometric immersion
$\tilde{f}\colon\,M^n\rightarrow \R^{n+2k_0-1}$
such that the subbundles  $N_s^{\tilde{f}}$ and $N_s^f$, $1\leq s\leq k_0$, endowed with the induced connection, correspond by a parallel isometry.
\end{proposition}

\proof Consider the normal  subbundle
$\mathcal{T}=N_1^f\oplus\ldots\oplus N^f_{k_0}$ with the induced connection
$
\hat{\nabla}^\perp_Y\eta=(\nabla^\perp_Y\eta)_{\mathcal{T}}.
$
We have to show that $\alpha_f$ still satisfies the Gauss, Codazzi and Ricci equations. In fact, the Gauss and
Codazzi equations are trivially satisfied. By Propositions~\ref{Lemma1} and \ref{proj2}, the subbundles
$\mathcal{T}$ and $\mathcal{T}^\perp$ are parallel in the normal connection along $Z$. Given $\eta\in
\mathcal{T}$, a simple calculation yields
$$
\hat{R}^\perp(X,Z)\eta-R^\perp(X,Z)\eta
=-\left(\nabla^\perp_X\nabla^\perp_Z\eta\right)_{\mathcal{T}^\perp}
+\nabla^\perp_Z\left(\nabla^\perp_X\eta\right)_{\mathcal{T}^\perp}
+\left(\nabla^\perp_{[X,Z]}\eta\right)_{\mathcal{T}^\perp}.
$$
Since  $R^\perp(X,Z)\eta \in \mathcal{T}$ by the Ricci equation, the left hand side vanishes and thus
$$
\hat{R}^\perp(X,Z)\eta=R^\perp(X,Z)\eta.
$$
Now using  Proposition \ref{nparalelos} we conclude that the Ricci equation is satisfied. Since $M^n$ is simply connected, the result follows from the Fundamental theorem of submanifolds. \qed\vspace{1,5ex}

Finally, we are in condition to prove Theorem \ref{main}.\vspace{1ex}

\proof  Assume that $f$ is nonregular. By Proposition \ref{proj2} there exists $k_0\leq \tau^f-1$ such that
$\xi_2^{k_0}=0$. By Proposition \ref{redcod},  there is a regular parabolic submanifold
$\tilde{f}\colon\,M^n\rightarrow \R^{n+2k_0-1}$ with $\xi_2^{\tau^{\tilde{f}}}=0$. It follows from Proposition
\ref{regrada} that $f$ is ruled. \qed

\section{Ruled parabolic}

The simple structure of ruled parabolic submanifolds allows us to give a parametric description of these
submanifolds. Using this description, we conclude that this submanifolds are generically regular. Then, we show that
ruled parabolic submanifolds are the only parabolic submanifolds that admit isometric immersions as
hypersurfaces.\vspace{1.5ex}

Let $v\colon\,I\subset\R\to\R^N$ be a smooth curve parametrized by arc length in some interval. Set  $e_1=dv/ds$ and
let $e_2,\ldots, e_{n-1}$ be  orthonormal normal vector fields along $v=v(s)$ parallel in the normal connection
of  $v$ in $\R^N$. Thus, 
\be\label{condi0} 
\frac{de_j}{ds}=b_je_1,\;\;\; 2\le j\le n-1, 
\ee 
where $b_j\in
C^\infty(I)$. Set $\Delta=\spa\{e_2,\ldots, e_{n-1}\}$ and let $\Delta^\perp$ be the orthogonal  complement in
the normal bundle. Take $e_0\in\Delta^\perp$ along $v$ such that
$$
P=\{e_0,(de_1/ds)_{\Delta^\perp}\}\subset\Delta^{\perp}
$$ 
satisfy that
\be\label{dim} 
\dim P = 2 
\ee 
and that $P$ is nowhere parallel in $\Delta^\perp$ along $v$, that is,
\be\label{dim2} 
\spa\{(de_0/ds)_{\Delta^\perp},  (d^2e_1/ds^2)_{\Delta^\perp}\} \not\subset P. 
\ee

We parametrize a ruled submanifold $M^n$ by
\be\label{para} 
f(s,t_1,\ldots,t_{n-1})={c(s)}+\sum_{j=1}^{n-1}t_j e_j(s)
\ee
where $(t_1,\ldots,t_{n-1})\in\R^{n-1}$ and $c(s)$ satisfies $dc/ds=e_0$.
To see that $f$ is parabolic, first observe that
$$
TM=\spa\{f_s\}\oplus\spa\{e_1\}\oplus\Delta
$$
where
$f_s=e_0+t_1de_1/ds+\sum_{j\ge 2}t_jb_je_1$.
Consider the orthogonal decomposition
\be\label{eum}
\left(\frac{de_1}{ds}\right)_{\Delta^\perp}= a_1e_0 + \eta.
\ee
Thus $\eta(s)\neq 0$ for all $s\in I$ from (\ref{dim}). Hence,
\be\label{tangente}
TM=\spa\{e_0+t_1(a_1e_0+\eta)\}\oplus\spa\{e_1\} \oplus\Delta.
\ee
Since $f_{st_j}=b_je_1\in TM,\; 2\le j\le n-1$,
we have that $\Delta\subset\Delta_f$. It follows easily from (\ref{eum}), (\ref{tangente}) and $\eta(s)\neq 0$ that
$$
f_{st_1}=\frac{de_1}{ds}\not\in TM.
$$
It is easy to see that $f_{ss}\not\in \spa\{f_{st_1}\}\oplus TM$, i.e.,  $\dim N_1^f=2$, is equivalent to
$$
\left(\frac{de_0}{ds}\right)_{\Delta^\perp} +t_1\left(\frac{d^2e_1}{ds^2}\right)_{\Delta^\perp} \not\in P.
$$
It follows that $\Delta=\Delta_f$. Therefore $f$ is parabolic in, at least, an open dense subset of $M^n$.

Let $f\colon M^n\to\R^N$ be a ruled parabolic submanifold and $\{e_2,\ldots,e_{n-1}\}$ an orthonormal frame for
$\Delta_f$ along an integral curve $c=c(s),\; s\in I$, of the unit vector field $X$ orthogonal to the rulings.
Without loss of generality (see Lemma 2.2 in \cite{BDJ}) we may assume that
$$
\frac{de_j}{ds}\perp\Delta_f,\;\;2\le j\le n-1.
$$
Now parametrize $f$ by (\ref{para}), where $e_0=X$ and $e_1=Z$.
That $f_{st_j}\in TM$ implies
\be\label{condi11}
\frac{de_j}{ds}\in\spa\{e_1,f_s\},\;\;\; 2\le j\le n-1.
\ee
Taking $t_1=0$, we obtain that
\be\label{condi22}
\frac{de_j}{ds}=a_je_0+b_je_1,\;\;\; 2\le j\le n-1,
\ee
where $a_j,b_j\in C^\infty(I)$.  Since  $\dim N^f_1=2$, we have
\be\label{uso}
\frac{de_1}{ds}= a_1e_0 + (de_1/ds)_{\Delta} + \eta
\ee
where $\eta\perp \spa\{e_0,e_1\}\oplus\Delta$ satisfies
$\eta(s)\neq 0$. Thus (\ref{condi11}) reduces to
$$
a_je_0\in\spa\{(1+t_1a_1+\ldots
+t_{n-1}a_{n-1})e_0+t_1\eta\} ,\;\;\; 2\le j\le n-1.
$$
Therefore $a_j=0$. From (\ref{condi22}) we have
$de_j/ds=b_je_1$ for $2\le j\le n-1$.\vspace{1ex}

We have proved the following result.

\begin{proposition}\po Let $c\colon\, I\subset \R\rightarrow \R^N,\;N-n\geq 2$, be a smooth curve. Let $\{e_0=dc/ds, e_1(s),\ldots,e_{n-1}(s)\}$ be orthonormal fields satisfying (\ref{condi0}), (\ref{dim}) and (\ref{dim2}) at any point. Then, the submanifold parametrized by
\be\label{paramregrada}
f(s,t_1,\ldots,t_{n-1})={c(s)}+\sum_{j\ge 1}t_j e_j(s)
\ee
where $(t_1,\ldots,t_{n-1})\in\R^{n-1}$, defines a ruled submanifold, that is  parabolic in an open dense subset of $M^n$.
Conversely, any ruled parabolic submanifold can be parametrized as in (\ref{paramregrada}).
\end{proposition}

 Let $f$ be a ruled parabolic submanifold  parametrized by (\ref{paramregrada}).  Assume that $f$ has critical index
$k-1=\tau_0^f$. The condition $\dim N_k^f=1$ is equivalent to \be\label{equiv1}
\frac{d^ke_1}{ds^k}\in
TM\oplus \spa\left\{\frac{d^{\ell-1}e_1}{ds^{\ell-1}}, \frac{d^{\ell-1} e_0}{ds^{\ell-1}}
+t_1\frac{d^\ell e_1}{ds^\ell},\;\; 2\le\ell\le k \right\}
\ee
where $TM$ was given by (\ref{tangente}). In particular, for
$t_1=0$ and using (\ref{uso}) we have
\be\label{equiv2}
\frac{d^{k-1}(a_1e_0 + \eta)}{ds^{k-1}}\in TM\oplus
\spa\left\{\frac{d^{\ell-2}(a_1e_0 + \eta)}{ds^{\ell-2}}, \frac{d^{\ell-1} e_0}{ds^{\ell-1}},\;\; 2\le\ell\le k\right\}
\ee
where now $TM=\spa\{e_0,e_1\}\oplus\Delta$.

It is easy to see that  (\ref{equiv1}) and (\ref{equiv2}) are   equivalent. In fact,  in (\ref{equiv2}) taking $\ell=2$ we obtain
that $\eta$ belongs to the subspace. If (\ref{equiv2}) is satisfied, it follows that the subspace in (\ref{equiv1}) is independent of the parameter $t_1$.  In particular, this shows again that
$\dim N_k^f=1$ is equivalent to $\xi_2^k=0$. Finally, we have that (\ref{equiv2}) is equivalent to
$$
\frac{d^{k-1}\eta}{ds^{k-1}}\in
\spa\left\{e_0,\frac{d e_0}{ds},
\ldots,\frac{d^{k-1} e_0}{ds^{k-1}},\eta,\ldots, \frac{d^{k-2}\eta}{ds^{k-2}}
\right\}\oplus \Delta.
$$
It is now clear that (\ref{equiv1}) will not be satisfied in general. In that sense and recalling Theorem \ref{main}, we can say that the  parabolic submanifolds are {\it generically\/} regular.

\begin{remark}\po
{\em A condition for a ruled regular parabolic submanifold in odd codimension to  satisfies  $\xi_2^{\tau^f}=0$ is the following:
$$
\frac{d^{\tau^f-1}\eta}{ds^{\tau^f-1}}\in
\spa\left\{e_0,\frac{d e_0}{ds},
\ldots,\frac{d^{\tau^f-2} e_0}{ds^{\tau^f-2}},\eta,\ldots, \frac{d^{\tau^f-2}\eta}{ds^{\tau^f-2}}
\right\}\oplus \Delta.
$$
}\end{remark}

Next we extend  the characterization of ruled parabolic submanifolds in codimension two given in \cite{DajFlo2} to arbitrary codimension. \vspace{1ex}

\begin{definition}\po{\em We say that a submanifold  $f\colon\,M^n\rightarrow \R^N$ is of \emph{surface type} if either $f(M)\subset L^2 \times \R^{n-2}$ where $L^2 \subset \R^{N-n+2}$ or \mbox{$f(M)\subset CL^2\times \R^{n-3}$} where $CL^2 \subset \R^{N-n+3}$ is a cone over a spherical surface $L^2\subset\mathbb{S}^{N-n+2}$.}
\end{definition}

\begin{theorem}\po\label{teo:hiper}
Let $f\colon\, M^n\rightarrow \R^N$ be a  ruled parabolic submanifold. If $M^n$ is simply connected then it admits an
isometric immersion as  a ruled hypersurface in $\R^{n+1}$ with the same rulings.
Conversely, if  $M^n$ admits an isometric immersion as a hypersurface in $\R^{n+1}$ and $f$ is not of surface type in any open subset, then $f$ is ruled.
\end{theorem}

\proof To prove the converse, assume that there exists an isometric immersion $g\colon\,M^n\rightarrow \R^{n+1}$ with Gauss map $N$.
We first show that
\be\label{same}
\Delta_g=\Delta_f.
\ee
Let $\beta\colon\,T_xM\times T_xM\rightarrow\R\<\eta_1\>\oplus\R\<N\>=\R^2$ be the symmetric bilinear form
$$
\beta(Y,V)=(\< A^f_{\eta_1}Y,V\>,\< A^g_N\,Y,V\>)
$$
where $\{\eta_1,\eta_2\}$ is as in (\ref{forma}). By the Gauss equation,  $\beta$ is flat with respect to the Lorentzian metric in $\R^2$ defined as
$\|\eta_1\|^2=1=-\|N\|^2$ and $\<\eta_1,N\>=0$, that is,
$$
\<\beta(X,Y),\beta(V,W)\>
-\<\beta(X,W),\beta(V,Y)\>=0.
$$
If (\ref{same}) is not  satisfied, and since $\dim\,\Delta_g\le n-2$, it follows easily that
$$
S(\beta)=\spa\{\beta(Y,V) : \;Y, V\in T_xM\}
$$
satisfies $S(\beta)=\R^2$.  From  Corollary 1 in \cite{Moore} we have $\dim N(\beta)=n-2$ where
$$
N(\beta)=\{Y\in T_xM : \beta(Y,V)=0,\;\;V\in T_xM\}.
$$
But since
$N(\beta)=\Delta_g\cap\Delta_f,$
it follows that (\ref{same}) holds.

Let
$$
A^g_N|_{\Delta^\perp}
=\left[\begin{array}{cc}
\bar{a} & \bar{b}\\
\bar{b} & \bar{c}
\end{array}\right].
$$
  From (\ref{cete}) we have
$$
C_T=\left[\begin{array}{cc}
m & 0 \\
n & m
\end{array}\right]
$$
for any  $T\in\Delta$. On the other hand,
$$
A^g_N \circ C_T=\left[\begin{array}{cc}
\bar{a}m+bn &\bar{b}m \\
\bar{b}m+\bar{c}n & \bar{c}m
\end{array}\right] .
$$
The symmetry of $A^g_N\circ C_T$ allows to conclude that $\bar{c}n=0$. Since $f$ is nowhere of surface type, it
follows from Lemma $6$ in \cite{DjFlToj}  that $n\neq 0$ for some $T\in\Delta$ in an open dense subset of $M^n$.
Thus $\bar{c}=0$ and therefore, by the Gauss equation, we may assume that $\bar{b}=b$.

The Codazzi equation for $A^f_{\eta_1}$ gives 
$$ 
\nabla_X bX -\<\nabla_X Z,X\>
(aX+bZ)-\nabla_Z(aX+bZ)+\<\nabla_ZX,Z\> bX+\< \nabla^\perp_Z \eta_1, \eta_2\> cX=0. 
$$
Taking the $Z$-component yields 
\be\label{codazzieta1} 
2b\<\nabla_X X,Z\>-a\< \nabla_Z X,Z\> -Z(b)=0. 
\ee 
The Codazzi equation for
$A^g_N$, that $\bar{c}=0$ and $\bar{b}=b$ give 
$$
\nabla_X bX-\< \nabla _X Z,X\>(\bar{a}X+bZ)
-\nabla_Z(\bar{a}X+bZ)+\< \nabla_ZX,Z\> bX=0. 
$$
Taking the $Z$-component yields 
\be\label{codazzig}
2b\<\nabla_X X,Z\>-\bar{a}\< \nabla_Z X,Z\>-Z(b)=0. 
\ee 
Subtracting  (\ref{codazzieta1}) from (\ref{codazzig}),
gives
$(a-\bar{a})\< \nabla_Z Z,X\>=0.$
If $\< \nabla_Z Z,X\>=0$, then $f$ is ruled. Thus, we may assume that
$a=\bar{a}$.  Now taking the $X$-component in both Codazzi equations
yields
$$
X(b)-a\<\nabla_XZ,X\>-Z(a)+2b\<\nabla_Z
X,Z\>+c\<\nabla^\perp_Z\eta_1,\eta_2\>=0
$$
and
$$
X(b)-a\<\nabla_XZ,X\>-Z(a)+2b\<\nabla_Z X,Z\>=0.
$$
It follows from the last two equations that
\be\label{eta2paralelo}
\<\nabla^\perp_Z\eta_1,\eta_2\>=0,
\ee
and we conclude from (\ref{normalregrada}) that $f$ is ruled. \vspace{1ex}

 We now prove the direct statement. In view of (\ref{forma}), we consider the tensor $A:TM\rightarrow TM$ where
$\Ker\, A=\Delta$ and
$$
A|_{\Delta^\perp}=\left[\begin{array}{cc}
  a & b \\
  b & 0
\end{array}\right].
$$
Since (\ref{eta2paralelo}) holds by assumption, it is easy to see that the tensor $A$ satisfies the  Gauss and Codazzi equations as a hypersurface, and this concludes the proof. \qed

\begin{corollary}\po\label{teo:reduzir1}
Let $f\colon\,M^n\rightarrow \R^{N}$ be a simply connected parabolic submanifold. Assume that there is $2\leq k_0\leq\tau^f-1$ such that $\dim N_{k_0}^f=1$. Then $f$ is ruled and $M^n$ admits an isometric immersion as a ruled hypersurface.
\end{corollary}

\proof We know from Proposition \ref{redcod} that there exists a regular parabolic isometric immersion
$\tilde{f}\colon\,M^n\rightarrow\R^{n+2k_0-1}$ such that $\xi_2^{k_0}=0$. It follows from Theorem~\ref{regrada}
that  $f$ is ruled. The result follows from Theorem \ref{teo:hiper}. \qed

\section{Nonruled parabolic submanifolds}

In this section we  study parabolic surfaces. First we show
that they are associated to parabolic differential equations.  
Then we give a complete characterization of their $s$-cross sections.

\medskip

Let $L^2$ be a Riemannian manifold endowed with a global system of coordinates. Then, let $f\colon\,L^2 \rightarrow \Q_\epsilon^N \subset \R^{N+\epsilon}$ where \mbox{$\epsilon=0,1$} and $N \geq 4$, be a surface of the sphere or the Euclidean space whose coordinate functions are linearly independent solutions  (of length 1 if  $\epsilon=1$) of the  parabolic equation
\be\label{PDE}
\frac{\partial^2 u}{\partial z^2}+W(u) +\epsilon \lambda u=0
\ee 
where $ W\in TL$ and $\lambda\in C^\infty(L^2)$.
If $\epsilon=0$, then  (\ref{PDE})  is equivalent to
$$
\tilde{\nabla}_Zf_\ast Z+f_\ast W=0
$$
where $Z=\partial/\partial z$. Thus $\alpha_f(Z,Z)=0$. If $\epsilon=1$, we have
$$
\tilde{\nabla}_Zf_\ast Z+f_\ast W+\lambda f=0
$$
and again $\alpha_f(Z,Z)=0$. In both situations $f$ is parabolic with  $Z$ asymptotic. 

Conversely, let $f\colon\,L^2 \rightarrow \Q_\epsilon^N$ be
parabolic endowed with the induced metric
and coordinates $(x,z)$ such that $\partial/\partial z=Z$ is asymptotic. The latter means that the
coordinate functions of $f$ satisfy (\ref{PDE}) with  $W=-\nabla_ZZ$ and $\lambda=\|Z\|^2$.
\medskip

Let $g\colon\,L^2\rightarrow \Q_\epsilon^N$ be a  parabolic surface and $\Sigma$ the vector space of classes of
functions $u\in C^\infty(L)$ that satisfy (\ref{PDE}), where for $\epsilon = 0$ we identify two  functions when
they differ by a  constant. Consider $L^2$ with the induced metric by $g$. Then (\ref{PDE}) takes the form
\be\label{PDE1} 
\hess_u(Z,Z)+\epsilon u=0 
\ee 
where $Z\in TL$ is an unit asymptotic field.

Given a parabolic submanifold $f:M^n\rightarrow \Q_\epsilon^N$, we denote
$$
\tau_*^f \;= \left\{\begin{array}{l}
\tau\;\;\;\;\;\;\;\;\;\;\; \mbox{if}\;\; N-n\;\;\; \mbox{is even}
\vspace{1.5ex}\\
\tau -1\;\;\;\;\; \mbox{if}\;\; N-n \;\;\; \mbox{is odd}.
\end{array} \right.
$$

Let $\Gamma_r$, $1\leq r\leq \tau_*^g$, be the vector space of classes of $r$--cross sections of $L^2$ where we
identify two sections when, up to a constant, they differ by a section of $N_{r+1}^g\oplus\ldots\oplus
N_{\tau^g}^g$.
Take $[\,h\,]\in\Gamma_r$ with $r<\tau_*^g$ and $1\leq r<s\leq \tau_*^g$. Then, set
$\mathcal{P}_r(h)=(\mu_1,\mu_2)\in\niu_r$. By  Corollary \ref{Proposition:isom}, there  exists an unique section
$\gamma_{r+1}\in N^g_{r+1}$ such that
$$
\mathcal{P}_r(h)=\mathcal{P}_r(-\gamma_{r+1}).
$$
Thus $\bar{h}_{r+1}=h+\gamma_{r+1}$ satisfies that
$\bar{h}_{r+1}=h+\gamma_{r+1} \in \Gamma_{r+1}$.
Using the above argument, it follows easily that there exist unique sections \mbox{$\gamma_j\in N_j^g $,} $r+1\leq j\leq
s$, such that 
\be\label{decomp} 
\bar{h}=h+\gamma_{r+1}+\ldots+\gamma_s 
\ee 
satisfies $\left[\,\bar{h}\,\right]\in\Gamma_s$.

We  show next that all the $\Gamma_r$'s  are isomorphic to $\Sigma$.
Given $\left[\,h \,\right]\in\Gamma_r$, set
$$
h=\epsilon\varphi g+ W+\delta
$$
where $W\in TL$, $\delta \in T^{\perp}L$ and $\varphi\in C^\infty(L)$ if $\epsilon=1$.
Given $Y\in TL$, we have
$$
h_\ast(Y)=\epsilon((Y(\varphi)-\< Y,W \>)g+\varphi Y)+\nabla_Y W+\alpha_g(Y,W)-A^g_\delta(Y)+\nabla^{\perp}_Y
\delta.
$$
Since the $TL$-component of $h_\ast(Y) $ vanishes, we obtain \be\label{pde0} 
\epsilon\varphi Y+\nabla_Y
W=A^g_\delta Y. 
\ee 
In particular, the map $(Y,U)\mapsto\<\nabla_Y W, U\> $ is symmetric. Thus, if $\epsilon=0$ and setting
$\Theta(U)=\< W,U \>$,
we have $d\Theta(Y,U)=0$.
Thus $W=\nabla\varphi,\;\textrm{for} \; \varphi \in C^\infty(L^2)$. If $\epsilon=1$, that the  $\spa
\{g\}$-component  of $h_{\ast}(Y)$ vanishes gives
$Y(\varphi)=\< Y,W \>,$
and again $W=\nabla\varphi$. In both cases, we obtain from (\ref{pde0}) we that
\be\label{sec}
\hess_\varphi+\epsilon\varphi I=A^g_\delta.
\ee

Consider the linear map $\Upsilon\colon\,\Gamma_r\rightarrow\Sigma$ defined by
$\Upsilon([h])=[\varphi]$.
Assume that $\Upsilon([h])=0$. Then $(h)_{T_gL}=\nabla\varphi=0$.
  From (\ref{sec}) we obtain $A^g_\delta=0$, which means $(h)_{N_1^g}=0$. Using $(iii)$ in Corollary
\ref{Proposition:isom} we obtain $h\in N_{r+1}^g\oplus\ldots\oplus N_{\tau^g}^g$. We conclude from the
definition of $\Gamma_r$ that $\Upsilon$ is injective.

Take $\varphi\in\Sigma$ and set 
$$
\mathcal{S}=\{\psi\in L_{sim}(TL,TL): \<\psi Z,Z\>=0\}. 
$$ 
Let $\Phi\colon\,N_1^g\rightarrow \mathcal{S}$ be the injective linear map defined by
$\Phi(\upsilon)=A^g_{\upsilon}$. From (\ref{PDE1}) and $\dim N^g_1=2$, we have that $\Phi$ is an isomorphism. It
follows that there exists a unique $\gamma_1\in N_1^g$ such that
$$
A_{\gamma_1}^g=\hess_\varphi+\epsilon \varphi I.
$$
We define
$\hat h=\epsilon \varphi g + \nabla\varphi +\gamma_1$.
Then,
$$
\hat h_*X=\epsilon X(\varphi)g+\epsilon\varphi X+\tilde{\nabla}_X \nabla\varphi + \tilde{\nabla}_X \gamma_1=\alpha_g(X,\nabla\varphi)+\nabla^\perp_X\gamma_1,
$$
and thus  $[\,\hat h\,]\in\Gamma_1$. We conclude from (\ref{decomp}) that $\Upsilon$ is an isomorphism. In this
way, we obtain the following recursive procedure for the construction of the $r$--cross sections for the parabolic
surfaces.

\begin{proposition}\po\label{proposition:rseccao} Let $g\colon\,L^2\rightarrow\Q_\epsilon^N$ be a regular parabolic surface. Then, any $r$--cross section, $1\leq r\leq\tau_*^g$ can be written as
\be\label{seccao} 
h_\varphi=\epsilon\varphi g+ g_*\nabla\varphi+\gamma_0+\gamma_1+\cdots+\gamma_r, 
\ee 
where
$\varphi$ satisfies (\ref{PDE}) and is unique (up to a constant if $\epsilon=0$), $\gamma_0$ is any section of
$N^g_{r+1}\oplus\ldots\oplus N_{\tau^g}^g$, $\gamma_1 \in N_1^g$ is the unique solution of
$A_{\gamma_1}^g=\hess_\varphi+\epsilon \varphi I$ and $\gamma_j \; , \; 2\leq j \leq r$, are the unique sections
given by (\ref{decomp}). Conversely, any function $h_\varphi$ with the form (\ref{seccao}) is a $r$--cross section
to $g$.
\end{proposition}

\section{The parametrizations}

In this section, we provide a parametrically description of all regular parabolic Euclidean submanifolds.  There are two alternative representation, the polar and bipolar parametrizations,
each of which is determined by a parabolic surface and a solution of a differential equation.
\medskip

Our  starting point, is to show how to construct parabolic submanifolds using   parabolic surface with non vanishing normal vector $\xi_2^\tau$, in particular,  any nonruled parabolic surface.
\medskip

Let  $g\colon\,L^2\rightarrow \Q_\epsilon^{N}$ a parabolic surface with $Z \in TL$ asymptotic and whose  normal
vector field $\xi_2^{\tau^g}$ does not vanish at any point.  Let $h$ be a $s$--cross section to g and  $\Lambda_s=N^g_{s+1} \oplus \ldots \oplus N^g_{\tau^g}$ for $1 \leq s \leq \tau_*^g $. Let
$\Psi\colon\,\Lambda_s \rightarrow\R^{N+\epsilon}$ be the map 
$$
\Psi(\delta)=h(x)+\delta
$$
where $\delta \in \Lambda_s(x)$.

\begin{proposition}\po\label{parabolicagerada}
At regular points, $M^n=\Psi(\Lambda_s)$ is a regular parabolic submanifold. Moreover,  $M^n$ is nonruled if $g$ is nonruled.
\end{proposition}

For the proof  we use the following general results.

\begin{lemma}\po\label{prop:cod1}
Let $f\colon\,M^n \rightarrow \R^N$
be a  parabolic submanifold. Then, we have:
\begin{itemize}
\item[(i)] If $\dim N^f_{k+1}=2$, then there exists
$\eta\in N^f_{k+1}$ such that the  components of
$\mathcal{P}_k(\eta)$ form a base of $N^f_k$.
\item[(ii)] Suppose that $N-n$  is odd,
$\dim N^f_{\tau^f-1}=2$ and that $\xi_2^{\tau^f}$ never vanishes. Then
$\mathcal{P}_{\tau^f-1}(\xi_2^{\tau^f})$ is a base of $N^f_{\tau^f-1}$.
\end{itemize}
\end{lemma}

\proof We prove $(i)$. From Corollary \ref{Proposition:isom} we have that $\mathcal{P}_k\vert_{N^f_{k+1}}$ is an
isomorphism and from Lemma \ref{Lemma:proj} that $\dim N_k^f=2$. Since $N^f_k$ has dimension $2$, there exists at
least one vector $(\mu_1,\mu_2)\in\niu_k$ with $\mu_2\neq 0$. Thus $\mu_1\:\textrm{and}\:
\mu_2$ are linearly independent and form a base of $N^f_k$.

 For the proof of $(ii)$ it is sufficient to show that \mbox{$(\nabla^\perp_Z\xi_2^{\tau^f})_{N^f_{\tau^f-1}}\!\!\neq 0$}. If the vector field vanishes, from the definition of $\niu_{\tau^f-1}$ we have 
$\< \nabla^\perp_X \xi_2^{\tau^f}, \xi_2^{\tau^f-1}\>=0.$
Thus $\xi_2^{\tau^f}= 0$ from Proposition \ref{Lemma1}, 
and this is a contradiction.\qed

\begin{lemma}\po\label{lemma:seccao1}
Let $\beta\colon\,M^n\rightarrow \R^{N+\epsilon}$ a $s$--cross section to $f$, $1\leq s\leq\tau^f$. Then,
$$
(\tilde\nabla_{Z}\beta_{\ast}(Z))_{N^f_{s-1}}=0.
$$
\end{lemma}

\proof For $s\geq 2$, we have that $\< \beta_{\ast}(Z),\xi_2^{s-1}\>=0.$  Then,
\begin{eqnarray*}
0\!\!&=&\!\!Z\< \beta_{\ast}(Z),\xi_2^{s-1}\>
=\<\tilde\nabla_{Z}\beta_{\ast}(Z),\xi_2^{s-1}\>+ \<
\beta_{\ast}(Z),\alpha^{s+1}(Z,Z,X,\ldots,X)\>\\
\!\!&=&\!\!\<\tilde\nabla_{Z}\beta_{\ast}(Z),\xi_2^{s-1}\>.
\end{eqnarray*}
Using Lemma \ref{Lemma:proj1}, is easy to prove  by a similar argument that
$$
\<\tilde\nabla_{Z}\beta_{\ast}(Z),\xi_1^{s-1}\>=0.
$$
For $s=1,$ since   $N_0^f=\Delta^\perp, \;\xi_1^0=X\; \mbox{and}\;\: \xi_2^0=Z$, the proof follows easily. \qed\vspace{1ex}

We now prove Proposition \ref{parabolicagerada}. \vspace{1ex}

\proof  Take a coordinate system $(x,z)$ of $L^2$ such that $Z=\partial/\partial z$ is asymptotic and  let $\{\eta_1,\ldots,\eta_k\}$ be an orthonormal frame of $\Lambda_s$. We parametrize $M^n$ by
$$
\Psi(x,z,t_1,\ldots,t_k)=h(x,z) +\sum_{j=1}^k t_j\eta_j(x,z)
$$
where $k=N-2s$ and $(t_1,\ldots,t_k)\in\R^k$.
From Lemma \ref{prop:cod1}, we have  $TM =\Lambda_{s-1}$ and $\Delta_{\delta}=\Lambda_s$. We claim that
$\Psi_\ast(Z)$ is asymptotic, that is,
$\tilde{\nabla}_Z \Psi_*(Z)\in TM$.
In view of  (\ref{comp}) it is sufficient to show for $\upsilon \in N^g_{s-1}$ that
$\<\tilde\nabla_{Z} \Psi_\ast(Z),\upsilon\>=0.$ 
Let $X=\partial/\partial x \in TL$.
We have that
\begin{eqnarray*}
\<\tilde\nabla_{Z} \Psi_\ast(Z),\xi_1^{s-1}\> \!\!&=&\!\!\<\tilde\nabla_Z h_\ast(Z),\xi_1^{s-1}\>+
\sum_{j=1}^{k}t_j\<\tilde{\nabla}_{Z}\tilde{\nabla}_{Z}\eta_j,
\alpha_\Psi^s(X,\ldots,X)\> \\ 
\!\!&=&\!\!\<\tilde\nabla_Z h_\ast(Z),\xi_1^{s-1}\>-\sum_{j=1}^{k}t_j\<\eta_j,
\tilde\nabla_{Z}\alpha_\Psi^{s+1}(Z,X,\ldots,X) \>\\
\!\!&=&\!\!\<\tilde\nabla_Z h_\ast(Z),\xi_1^{s-1}\>.
\end{eqnarray*}
By a similar argument, we obtain
$$
\<\tilde\nabla_{Z} \Psi_\ast(Z),\xi_2^{s-1}\>=\<\tilde\nabla_Z h_\ast(Z),\xi_2^{s-1}\>.
$$
Now Lemma \ref{lemma:seccao1} and $N_{s-1}^g=\mbox{span}\{\xi_1^{s-1},\xi_2^{s-1}\}$ give the claim. Observe that it follows from  Lemma \ref{prop:cod1} that $N_k^\Psi =N^g_{s-k}$. This concludes the first part of the proof.

Assume that $g$ is nonruled. From Lemma \ref{Lemma:proj1}  we have that $\xi_2^s$ and $\Psi_\ast(Z)$ are
orthogonal. Being $\eta_s \in \Delta_{\Psi}^{\perp}=N_s^g$  a unit asymptotic vector field to $\Psi$, we obtain
that $\Psi$ is ruled if and only if $(\tilde{\nabla}_Z\eta_s)_{N^g_s}=0$. Now the proof follows from Corollary~\ref{coro:regrada}. \qed\vspace{1,5ex}

Our goal now is to show that any parabolic submanifolds with non vanishing normal vector field $\xi_2^\tau$,
in particular, all nonruled regular parabolic submanifolds,  can be locally parametrized by a parabolic surface
using Proposition \ref{parabolicagerada}.
\medskip

Given a parabolic submanifold $f\colon\,M^n\rightarrow\Q_{\epsilon}^N$, due to the  local nature  of our work,  we may
assume that $f$ is the saturation of a fixed cross section $L^2\subset M^n$ to the relative nullity foliation.
  From Proposition \ref{nparalelos}, each  $N^f_k$ can be viewed as a plane bundle along~$L^2$.

\begin{definition}\po{\em
Let $f\colon\,M^n\rightarrow\Q_{\epsilon}^{N-\epsilon}$ be a regular parabolic submanifold. A \emph{polar
surface} to $f$ is an immersion of a cross section  $L^2$  as above, defined as follows:
\begin{itemize}
\item[(i)]If $N-n-\epsilon$ is odd, then
$g\colon\,L^2\rightarrow\mathbb{S}^{N-1}$ is defined by
$$
\spa \{g(x)\}=N_{\tau^f}^f(x).
$$
\item[(ii)]If $N-n-\epsilon$ is even, then  $g\colon\,L^2\rightarrow\R^N$
is any surface such that
$$
T_{g(x)}L=N_{\tau^f}^f(x),
$$
up to parallel identification in $\R^n.$
\end{itemize}}
\end{definition}

\begin{proposition}\po\label{polar} Any regular parabolic submanifold  $f\colon\,M^n\rightarrow \Q_{\epsilon}^N$ with non vanishing normal vector field $\xi_2^{\tau^f}$ admits a polar surface $g$ locally. Moreover, $g$ is parabolic and nonruled if $f$ is nonruled and has no Euclidean factor.
\end{proposition}

We will use the following fact.

\begin{lemma}\po\label{vectparalelo}
Assume that $f$ has even codimension. Let $\eta\in N^f_{\tau^f}$ and
$$
\mu_1=(\tilde{\nabla}_X\eta)_{N^f_{\tau^f-1}},\;\;\; \mu_2=(\tilde{\nabla}_Z\eta)_{N^f_{\tau^f-1}}
$$ 
be such that
$\mu_2\neq 0$. Then,
$$
\niu_{\tau^f-1}=\{(a\mu_1+b\mu_2,a\mu_2):\; a,b \in C^\infty(M)\}.
$$
\end{lemma}

\proof Since $\<(\tilde{\nabla}_Z\eta)_{N^f_{\tau^f-1}},
\xi_2^{\tau^f-1}\>=\<\eta,\tilde{\nabla}_Z\xi_2^{\tau^f-1}\>=0$,  
the definition of $\niu_{\tau^f-1}$ and 
Lemma \ref{Lemma1} yield $(\mu_2,0)\in\niu_{\tau^f-1}$. Since 
$\dim N^f_{\tau^f-1}=2$, we easily 
conclude that $N^f_{\tau^f-1}=\mbox{span}\{(\mu_1, \mu_2)$, $(\mu_2,0)\}$, and the proof follows.\qed

\begin{remark}\po {\em Notice that  $\eta=\xi_2^{\tau^f}/\|\xi_2^{\tau^f}\|\in N_{\tau^f}^f$
 satisfies $(\tilde{\nabla}_Z\eta)_{N_{\tau^f-1}^f}\neq 0$. In fact, from Proposition \ref{Lemma1} it is easy to see that 
$\< \tilde{\nabla}_Z\eta, \xi_1^{{\tau^f}-1}\> \neq 0.$}
\end{remark}

We now prove Proposition \ref{polar}.\medskip

\proof In the case of odd codimension, the existence of a polar surface follows from $(ii)$ of Lemma
\ref{prop:cod1}. Assume that $\dim N_{\tau^f}^f=2$. Let $\{\eta_1,\eta_2\}$ be a  base of  $N_{\tau^f}^f$ constant
along $\Delta$. We show that there exist linearly independent $1-$forms, $\theta_1, \theta_2$ so that the
differential equation 
\be\label{eqdf} 
dg=\theta_1 \eta_1+\theta_2\eta_2 
\ee 
has solution.

Take a non vanishing asymptotic vector field $Z\in TM$  and consider the isomorphism
$P\colon\,\Delta^{\perp}\rightarrow TL.$ Let $U=P(Z)\in TL$ and $(u,w)$ a coordinate system  on $L^2$  such that
$U=\partial/\partial u.$ Set $W=\partial/\partial w\in TL$ and $X=P^{-1}(W) \in \Delta^\perp$. Endow $L^2$
with the metric which makes the base $\{U,W\}$ orthonormal and positively oriented.
Let $\eta_1, \eta_2 \in N_{\tau^f}$ be linearly independents vector fields constant along~$\Delta$. Without loss of generality, we my assume 
$\mu_2=(\tilde{\nabla}_Z\eta_1)_{N^f_{\tau^f-1}}\neq 0$.
According to Lemma \ref{vectparalelo}, there are $a,b\in C^\infty(M)$ with $b\neq 0$ such that 
\be\label{projeta}
\mathcal{P}_{\tau^f-1}(\eta_1)=(\mu_1,\mu_2)\quad\mbox{and}\quad
\mathcal{P}_{\tau^f-1}(\eta_2)=(a\mu_1+b\mu_2,a\mu_2). 
\ee 
Consider 1-forms
\be\label{dif} 
\theta_1=a^1du+a^2dw
\quad\mbox{e}\quad \theta_2=b^1du+b^2dw, 
\ee   
where $a^1,a^2,b^1,b^2 \in C^\infty(L^2)$. We show that
we can choose \mbox{$a^1,a^2,b^1,b^2 \in C^\infty(L)$} such that (\ref{eqdf}) has  solution . The integrability
condition for (\ref{eqdf}) is \bea 0\!\!&=&\!\! d\theta_1\eta_1+d\theta_2\eta_2
+\theta_1 \land d\eta_1+\theta_2\land d\eta_2\\
\!\!&=&\!\!d\theta_1\eta_1+d\theta_2\eta_2
+(a^1\frac{\partial \eta_1}{\partial w}
-a^2\frac{\partial\eta_1}{\partial du})dV
+(b^1\frac{\partial \eta_2}{\partial v}
-b^2\frac{\partial \eta_2}{\partial u})dV\\
\!\!&=&\!\!d\theta_1\eta_1+d\theta_2\eta_2 +(\tilde{\nabla}_{a^1W-a^2U}\:\eta_1
+\tilde{\nabla}_{b^1W-b^2U}\:\eta_2)dV 
\eea 
where  $dV$ stands for the volume element of $L^2$. Then, we must have

$$
(\tilde{\nabla}_{a^1W-a^2U}\:\eta_1
+\tilde{\nabla}_{b^1W-b^2U}\:\eta_2)_{N_{\tau^f-1}}=0.
$$

From (\ref{projeta}) we may rewrite the above equation as
\be\label{sistemadif1} \left\{
\begin{array}{l}
a^1+ab^1=0\vspace*{1ex}\\
a^2-bb^1+ab^2=0.
\end{array}
\right.
\ee
Then, let $e,\ell\in C^\infty(L)$ be such that
$$
\tilde{\nabla}_{a^1W-a^2U}\:\eta_1+\tilde{\nabla}_{b^1W-b^2U}\:\eta_2=e\eta_1+\ell\eta_2.
$$
We claim that there exist $a^1,a^2,b^1,b^2\in C^\infty(L)$ such that $\theta_1,\theta_2$ satisfy(\ref{sistemadif1}) and
$$
\left\{
\begin{array}{lll}
d\theta_1 \!\!&=&\!\!e\:dV\vspace*{1ex} \\
d\theta_2 \!\!&=&\!\!\ell\:dV,
\end{array}\right.
$$
or equivalently, 
\be\label{sistemadif3} 
\left\{
\begin{array}{lll}
a^2_u-a^1_w \!\!&=&\!\!e\vspace*{1ex} \\
b^2_u-b^1_w \!\!&=&\!\!\ell.
\end{array}\right.
\ee 
From (\ref{sistemadif1}) and (\ref{sistemadif3}) we have
$$
\left\{\begin{array}{l}
a^1=-ab^1\vspace*{1ex}\\
a^2=bb^1-ab^2\vspace*{1ex}\\
b_ub^1+bb^1_u-a_ub^2-a(b^2_u-b^1_w)+a_wb^1=e\vspace*{1ex}\\
b^2_u-b^1_w=\ell.
\end{array}\right.
$$
The  two last equations give
\be\label{sistemadif4} 
\left\{\begin{array}{l}
b_ub^1+bb^1_u-a_ub^2+a_wb^1=e+a\ell\vspace*{1ex}\\
b^2_u-b^1_w =\ell.
\end{array}\right.
\ee 
We assume  $a_u\neq 0$ without loss of generality. The first equation of
(\ref{sistemadif4}) yields
$$ 
b^2 =-\frac{1}{a_u}(e+a\ell-(b_u+a_w)b^1+bb^1_u). 
$$ 
We take  $b^1$ to be a solutions of the above linear parabolic equation (see p.\ 367 of \cite{Ev}), and now the claim follows easily.

If $f$ has a Euclidean factor, take $T$ a parallel subbundle of the relative nullity subbundle of $f$. It is easy
to see that under these conditions the subbundle $T\oplus \nabla^\perp\oplus N_1^g$ is a normal parallel subbundle
of $g$. Thus, the codimension of $g$ can be reduced. The converse is similar.

We claim that $g$ has an asymptotic vector. First observe that $N_1^g=N_{\tau_*^f-1}^f.$ Thus, in odd codimension,
we have from (\ref{eqdf}) and (\ref{sistemadif1}) that \be\label{derivadapolar} 
g_*\partial/\partial u
=a^1\eta_1+b^1\eta_2=-ab^1\eta_1+b^1\eta_2. 
\ee 
Therefore, in view of (\ref{projeta}) we obtain
$$
(\tilde{\nabla}_Z g_*\partial/\partial u)_{N^f_{\tau^f-1}} = -ab^1\mu_2+ab^1\mu_2=0.
$$
For even codimension, the claim follow from Lemma \ref{lemma:seccao1}. Hence $g$ is parabolic.

To complete the proof suppose that  $f$ is nonruled. 
We show that $g$ is also nonruled. If the codimension of f is odd,  since $\xi_2^{\tau^f}\neq 0$, then $TL$ is spanned by $\{(\tilde{\nabla}_X\xi_2^{\tau^f})_{N_{\tau_*^f}^f},
(\tilde{\nabla}_Z\xi_2^{\tau^f})_{N_{\tau_*^f}^f}\}$, being $(\tilde{\nabla}_Z\xi_2^{\tau^f})_{N_{\tau_*^f}^f}$
an asymptotic field. 

The definition of $\niu_{\tau_*^f}$ allows us to conclude that the unit asymptotic field $\gamma$ is normal to $\xi_2^{\tau_*^f}$
Then, $g$ is ruled if and only if
$(\tilde{\nabla}_Z \gamma)_{N_{\tau_*^f}^f}=0$.
Thus $g$ is nonruled by Corollary~\ref{coro:regrada}. 
In the even codimension case,  we have  
$$
N_1^g=\spa\{(\tilde{\nabla}_Z\eta_1)_{N_{\tau^f-1}^f},
(\tilde{\nabla}_X\eta_1)_{N_{\tau^f-1}^f}\}.
$$
 From (\ref{projeta}) and (\ref{derivadapolar}) it is easy to
conclude that 
\be\label{noregrada2} 
\xi_2^{1\:g}=b\mu_2 =b(\tilde{\nabla}_Z\eta_1)_{N_{\tau^f-1}}. 
\ee
 Let $\lambda=\|b(\tilde{\nabla}_Z\eta_1)_{N_{\tau^f-1}}\|^{-1}.$ 
It follows from  (\ref{normalregrada}) that $g$ is ruled if and only if  
$$
(\tilde{\nabla}_U\lambda (\tilde{\nabla}_Z\eta_1)_{N_{\tau^f-1}})_{N^f_{\tau^f-1}}=0.
$$
 From our assumption that
$\eta_1$ is constant along $\Delta_f$, it follows that
$$
0=(\tilde{\nabla}_U\lambda (\tilde{\nabla}_Z\eta_1)_{N_{\tau^f-1}^f})_{N_{\tau^f-1}^f}
=U(\lambda)(\tilde{\nabla}_Z\eta_1)_{N_{\tau^f-1}^f} +\lambda(\tilde{\nabla}_Z
(\tilde{\nabla}_Z\eta_1)_{N_{\tau^f-1}^f})_{N_{\tau^f-1}^f}.
$$
Thus,
$$
(\tilde{\nabla}_Z (\tilde{\nabla}_Z\eta_1)_{N_{\tau^f-1}}^f)_{N_{\tau^f-1}^f}
\in(\tilde{\nabla}_Z\eta_1)_{N_{\tau^f-1}^f}.
$$
Since $(\tilde{\nabla}_Z\eta_1)_{N_{\tau^f-1}}$ is normal to $\xi_2^{{\tau^f-1}\:f}$, we obtain 
$$
(\tilde{\nabla}_Z\; \xi_2^{{\tau^f-1}\:f}/ \| \xi_2^{{\tau^f-1}\:f}\|)_{N_{\tau^f-1}^f}=0,
$$
and conclude from Corollary \ref{coro:regrada} that $f$ is ruled. 
This is a contradiction.\qed\vspace{1,5ex}

The following is the \emph{polar parametrization}.

\begin{theorem}\po\label{theorem:polar} 
Given a parabolic surface $g\colon\,L^2\rightarrow \Q_\epsilon^N$ with non vanishing
normal vector $\xi_2^{\tau^g}$ and $1 \leq s \leq \tau_*^g$, consider the smooth map $\Psi\colon\,\Lambda_s
\rightarrow \R^N$ defined by
\be\label{paramet1} 
\Psi(\delta)=h(x)+ \delta 
\ee
where $\delta \in\Lambda_s=N_{s+1}^g\oplus\ldots\oplus N_{\tau^g}^g$ and $h$ is any $s$--cross section to $g$. Then, at
regular points,  $M^n=\Psi(\Lambda_s)$ is a regular parabolic submanifold with polar surface~$g$. Moreover, if
$g$ is nonruled, then
$M^n=\Psi(\Lambda_s)$ is nonruled.

Conversely, any parabolic submanifold $f\colon\,M^n\rightarrow\R^N$ without local Euclidean factor and with non vanishing normal vector
$\xi_2^{\tau^f}$ admits a local parametrization (\ref{paramet1}), where $g$ is
a polar surface to $f$.
\end{theorem}

\proof The direct statement follows from Proposition \ref{parabolicagerada}. For the converse, take a polar
surface $g\colon\,L^2\rightarrow \Q_\epsilon^N$ to $f$. It is easy to see that under these conditions that
$\Delta_f=\Lambda_{\tau_*^f}$ and $TM=\Lambda_{\tau_*^f-1}$ along $L^2$. Thus, the section $h=f_{\mid _{L^2}}$ is
a $\tau_*^f$--cross section to $g$.\qed\vspace{1,5ex}

Observe that picking a different $\gamma_0$ in (\ref{seccao}) only results in a
reparametrization of $\Psi(\Lambda_s)$. Hence, it is convenient to take $\gamma_0=0$ when using the recursive
procedure to generate $s$--cross sections. 

The polar parametrization is very effective for
submanifolds in low codimension since the recursive procedure has few iterations. For instance, in codimension
two it suffices to take a $1$--cross section of the form $h_\varphi=\nabla\varphi+\gamma_1,$ where $\gamma_1 \in
N^f_1$  is unique satisfying $A_{\gamma_1}=\hess_\varphi$ for a 
given solution $\varphi$ of (\ref{PDE}).

\begin{definition}\po {\em We define the \emph{bipolar surface} to a parabolic submanifold  $f$ to be any polar
surface to a polar surface to~$f$}.
\end{definition}
\begin{proposition}\po\label{bipolar1} Any nonruled parabolic submanifolds admits locally a bipolar surface.
\end{proposition}
\proof From Proposition \ref{polar}, $f$ admits locally a  nonruled polar surface $g$. Then,  Proposition
\ref{regrada} gives $\xi_2^{\tau^g}\neq 0$. The proof now follows from Proposition
\ref{polar}\qed

\begin{definition}\po {\em Let $g\colon\,L^2\rightarrow \Q_\epsilon^N$ be a parabolic surface and $0\leq s\leq\tau_*^g-1$.
 We call  \emph{dual s--cross section} to $g$ any element $h\in C^\infty(L^2,\R^{N+\epsilon})$ satisfying
$$
h_*(TL) \subset \epsilon\:\spa\{g\} \oplus N_0^g \oplus \ldots \oplus N_s^g
$$
at any point.}
\end{definition}

Notice that a dual $0$-section to a parabolic surface in Euclidean space is just a bipolar surface.

\begin{proposition}\po\label{dualsection} Let $g\colon\,L^2\rightarrow \Q_\epsilon^N$ be a regular parabolic surface with polar surface $\hat{g}$. Any dual $s$-section to $g$ is a $([N/2]-s-1)$-section to $\hat{g}$.
\end{proposition}

\proof We have $\tau_\ast^g=\tau_\ast^{\hat{g}}=[N/2]-1$ and
$N_s^g=N_{\tau_\ast^{\hat{g}}-s}^{\hat{g}}$. The proof follows easily. \qed\vspace{1,5ex}

The following is the \emph{bipolar parametrization}.

\begin{theorem}\po\label{bipolar} Given a parabolic surface $g\colon\,L^2\rightarrow \Q_\epsilon^N$ with non vanishing normal vector $\xi_2^{\tau^g}$
and $0\leq s\leq \tau_*^g-1,$ consider the smooth map  $\tilde{\Psi}\colon\,\tilde{\Lambda}_s \rightarrow \R^N$
defined by 
\be\label{dualparamet} 
\tilde{\Psi}( \tilde{\delta})=\tilde{h}(x)+\tilde{\delta}
\ee 
where $\tilde{\delta}\in\tilde{\Lambda}_s=\epsilon\:\spa\{g\} \oplus N_0^g \oplus \ldots
\oplus N_{s-1}^g$ and $\tilde{h}$ is any dual $s$--cross section to $g$.  Then, at regular points,
$M^n=\tilde{\Psi}(\tilde{\Lambda}_s)$ is a  nonruled parabolic submanifold with bipolar surface $g$.

  Conversely, any nonruled parabolic submanifold  $f\colon\,M^n\rightarrow\R^N $  without local Euclidean factor admits a local parametrization (\ref{dualparamet}), where $g$ is a bipolar surface to $f$.
\end{theorem}

\proof The result follows from  Theorem \ref{theorem:polar} and Propositions \ref{bipolar1} and \ref{dualsection}.
\qed\vspace{1,5ex}

Next, we give a simple way to parametrize parabolic submanifolds. \vspace{1ex}

Let $g\colon\,L^2\rightarrow \Q_\epsilon^N$ be a simply connected nonruled parabolic surface endowed with the
metric induced by $g$ and $\{X,Z\}$ an orthonormal tangent frame with $Z$  asymptotic. Let $J\in
End\,(TL)$ be defined by
$$
J(X)=Z\;\;\mbox{and}\;\; J(Z)=0
$$
and let $R\in End\,(TL)$ the reflection defined by $$ R(X)=X\;\;\mbox{and}\;\;R(Z)=-Z.
$$
Now consider the linear second order parabolic operator 
$$
L(\varphi)=ZZ(\varphi)+\Gamma_2
X(\varphi) -\Gamma_1 Z(\varphi)+(X(\Gamma_2) -Z(\Gamma_1)+(\Gamma_1)^2-(\Gamma_2)^2-\epsilon)\varphi 
$$ 
where $Y=[X,Z]=\Gamma_2 Z-\Gamma_1 X$. Let $\varphi\in C^\infty(L)$ satisfy $L(\varphi)=0$  and  let $\psi$ be the
$1$-form such that
$d\psi(X,Z)=-\varphi$.
\begin{lemma}\po The differential equation
\be\label{1equation} 
d\theta=d\varphi\circ J +\varphi Y^*\circ R+\epsilon\psi 
\ee 
is integrable.
\end{lemma}

\proof From our assumptions,  we easily obtain
$d^2\theta(X,Z)=-L(\varphi)$, and this concludes the proof. \qed

\begin{lemma}\po\label{lemabip} The differential equation
\be
\label{edpdual1}
 dh=\epsilon \psi g + dg\circ (\theta I+\varphi J )
\ee 
is integrable, where $\theta$ is a solution of (\ref{1equation}).
\end{lemma}

\proof An easy computation yields
\bea
d^2h(X,Z)\!\!\!&=&\!\!\!
\epsilon(d\psi(X,Z)+\varphi)g
+(d\theta(X)+\varphi\Gamma_1-Z(\varphi) -\epsilon\psi(X))Z\\
\!\!\!&-&\!\!\!(d\theta(Z)+\varphi\Gamma_2-\epsilon\psi(Z))X.
\eea
Thus, we conclude that $d^2h=0$. \qed

\begin{theorem}\po Let  $g\colon\,L^2\rightarrow
\Q_\epsilon^{N-\epsilon}$ a simply connected nonruled parabolic surface, $\varphi \in C^\infty(L)$ so that
$L(\varphi)=0$ and $h\colon\,L^2\rightarrow \R^N$ a solution of (\ref{edpdual1}). Then, the map
$\Psi\colon\,L^2\times \R^{2s-\epsilon}\rightarrow\R^N$
defined by,
$$
\Psi(x,t)=h(x)+\epsilon\, t_0g(x)+
\sum\limits_{j=1}^{s}\left(t_{2j-1}
\frac{\partial^jg}{\partial v\partial u^{j-1}}+
 t_{2j}\frac{\partial^jg}{\partial u^{j}} \right)(x)
$$
where $0\leq s\leq [(N-\epsilon)/2]-2$ and $(u,v)$ is a coordinate system of $L^2$ such that $\partial/\partial v$
is asymptotic, parametrizes, at regular points, a parabolic submanifold.

Conversely, any nonruled parabolic submanifold  without local Euclidean factor can be locally parametrized in this
way.
\end{theorem}

\proof It is clear for $0\leq j \leq \tau_*^g$ that 
$$
N_j^g=\spa\left\{\left(\frac{\partial^{j+1} g}{\partial u^j \partial v}\right)_{N_j^g},\left(\frac{\partial^{j+1}g}
{\partial u^{j+1} }\right)_{N_j^g}\right\}. 
$$
In (\ref{dualparamet}) we take  $\tilde{h}$ to be a dual
$0$--cross section to $g$  without loss of generality. It remains to show that any dual $0$-section to $g$ can be written as a solution of
(\ref{edpdual1}).

Given a dual $0$-section $\tilde{h}$ to $g$,  we need a $1$-form $\Psi$ and $S\in End\,(TL)$ such~that
$$
d\tilde{h}=\epsilon \Psi g+dg\circ S.
$$
An easy computation yields
\begin{eqnarray*}
d^2\tilde{h}(X,Z)
\!\!&=&\!\!\epsilon (d\psi(X,Z)-\< X,SZ\>+\< Z,SX\>)g+(\nabla_X S)Z +\alpha_g(X,SZ)\\
\!\!&&\!\!-(\nabla_Z S)X-\alpha_g(Z,SX)+\epsilon(\psi(Z)X-\psi(X)Z).
\end{eqnarray*}
Thus, the integrability conditions reduces to the equations \be\label{dual1} 
\alpha_g(X,SZ)=\alpha_g(Z,SX), 
\ee
\be\label{dual2} 
(\nabla_X S)Z-(\nabla_Z S)X=\epsilon(\psi(X)Z-\psi(Z)X), 
\ee 
and for $\epsilon =1$ the additional equation 
\be\label{dual3} 
d\psi(X,Z)=\< SZ,X\>- \< SX,Z\>. 
\ee 
From (\ref{dual1}) and since
$\alpha_g(X,X)$ and $\alpha_g(X,Z)$ are linearly independent,  
we have 
$$
S=\theta I+ \varphi J
$$
where $\theta, \varphi \in C^\infty(L).$ The  left side of (\ref{dual2})  gives us 
$$ 
\nabla_X\theta
Z-\nabla_Z (\theta X+\varphi Z)
+\Gamma_1SX-\Gamma_2SZ\!=\!(d\theta(X)+\varphi\Gamma_1-d\varphi(Z))Z -(d\theta(Z)+\varphi\Gamma_2)X. 
$$
Thus (\ref{dual2}) is
equivalent to
\begin{displaymath}
\left\{
\begin{array}{rcl}
d\theta(X) \!\!&=&\!\!-\Gamma_1\varphi+d\varphi(Z)+\epsilon\psi(X)\vspace*{1ex}\\
d\theta(Z) \!\!&=&\!\! \< Y, -Z\>\varphi+\epsilon\psi(Z).
\end{array}
\right. 
\end{displaymath} 
Hence,
$$
d\theta=d\varphi\circ J +\varphi Y^\ast\circ R
+\epsilon\psi,
$$
and from (\ref{dual3}) we easily get
$d\psi(X,Z)=-\varphi$.
The result follows from Theorem~\ref{bipolar} and Lemma \ref{lemabip}. \qed

\section{The singularities}

In this section we show that the nowhere nonruled complete parabolic submanifolds  are  surface-like, that is,
they are isometric to  $L^2\times \R^{n-2}$. We also describe the singular set of nonruled parabolic submanifolds of dimension at least four.

\vspace{1ex}

The complete submanifolds  $f\colon\,M^n \rightarrow \R^N$ with rank $\rho\leq 2$, had been studied in
\cite{DajGro}. If $M^n$ does not contain an open set $L^3\times\R^{n-3}$ with $L^3$ unbounded, then the
following holds in the open set $M^*\subset M^n$ where $\rho=2$.
\begin{itemize}
\item[(i)] $M^*$ is an union of smoothly ruled strips.
\item[(ii)] If $f$ is completely ruled on $M^*$, then it is completely ruled everywhere and a cylinder on each component of
the complement of the closure of $M^*$.
\end{itemize}

A ruled submanifold is called \emph{completely ruled} if each leaf is a complete affine space. The leaves in
each connected component of $ M^n$,  called a \emph{ruled strip}, form an affine vector bundle over a curve
with or without end point \cite{DajGro}.\vspace{1ex}

Given a ruled parabolic submanifold $f\colon\,M^n \rightarrow \R^N$, let $\tilde M^n$ be the extension of
$f(M^n)$ (with possible singularities) obtained by extending each leaf to a complete affine Euclidean space
$\R^{n-1}$. We have the following result.

\begin{proposition}\po\label{regradacompleta}
Let $f\colon\,M^n \rightarrow \R^N$ a ruled parabolic submanifold. Then $\tilde M^n$  is a ruled strip. Moreover,
if  $c$ is complete and the function  $a_1$ defined in (\ref{eum}) satisfy $|a_1(s)|\leq K< +\infty$, then
$\tilde{M}^n$ is complete.
\end{proposition}

\proof Using (\ref{paramregrada}) we parametrize $M^n$ by 
$$ 
f(s,t_1,\ldots,t_{n-1})={c(s)}+\sum_{j\ge 1}t_j e_j(s)w
$$
where
$$ 
\frac{de_1}{ds}= a_1e_0 + \delta + \eta\;\;\;\mbox{and}\;\;\; \frac{de_j}{ds}=b_j e_1,\;\;2\le j\le n-1,
$$
$\delta=(de_1/ds)_{\Delta}$ and $\eta\perp \spa\{e_0,e_1\}\oplus\Delta$ is nonsingular  for every  $s\in I$.
We have,
$$
TM=\spa\{(1+t_1a_1)e_0+t_1\eta\}\oplus\spa\{e_1,\ldots,e_{n-1}\},
$$
and is now easy to conclude that $f$ is nonsingular. Thus $\tilde{M}^n$ is a ruled strip.

Next, suppose that $c$ is complete. Notice that
$$
\|f_s\|^2\ge (1+t_1a_1(s))^2+t_1^2\|\eta(s)\|^2.
$$
We claim that $\tilde M^n$ is complete. If $|t_1|\leq M<\infty$, from our assumption that $|a_1(s)|\leq K< \infty$ we obtain
$\|f_s\|^2\ge L>0$. On the other hand, it is easy to see that any divergent curve
$\gamma(u)=f(s(u),t_1(u),...,t_{n-1}(u)), u \in[0,+\infty)$, in $\tilde M^n$ with at least one $t_i, 1\leq i \leq
n-1$, unbounded has infinity length. Thus, any divergent curves in $\tilde M^n$ has infinity length, and the proof
follows. \qed\vspace{1.5ex}

Observe that any ruled parabolic submanifold parametrized by (\ref{paramregrada}) with \mbox{$b_j=0,\; 2\leq j\leq n-1$},
everywhere is a product $L^2\times \R^{n-2}$. On the other hand, if there exist $j\in \{2,\ldots, n-1\}$ such that
 $b_j\neq0 $ everywhere then the submanifold does not contain an open set $L^2\times \R^{n-2}$.

\begin{theorem}\label{maincompleta}\po
Let $f\colon\,M^n\rightarrow\R^N$, $n\ge 3$, be a complete submanifold which is nonruled in any open set and
parabolic in an open dense set $\mathcal{O}$. Then, any connected component of $\mathcal{O}$  is isometric to
$L^2\times \R^{n-2}$ and $f$ splits accordingly.
\end{theorem}

\proof From Lemma $6$ in \cite{DajGro} it is easy to see that
either $C=0$ or
\be\label{formaCT}
C_T=\left[\begin{array}{cc}
  0 & 0 \\
  n & 0
\end{array}\right]
\ee 
where $T\perp \mbox{Ker}\; C$. We have a disjoint decomposition
$\mathcal{O}=M_0\:\cup\:M_1$,
where  $M_0$ is the closet set where $C=0$. We now argue that the open set $M_1$ is empty. It follows from Lemma
$1.8$ in \cite{DajGro} that $M_0$ and $M_1$ are saturated, i.e. they are unions of complete leaves of $\Delta.$
We have from Lemma 1.5 in \cite{DajGro} and~(\ref{formaCT}) that
$$
0=(\nabla_XC_{T})Z-(\nabla_ZC_{T})X
=n\< \nabla_XZ,X\> Z
-Z(n)Z-n\<\nabla_ZZ,X\> X
$$
where $T\perp \ker C$ is an unit field. Therefore $\<\nabla_ZZ,X\> =0$, i.e., $M_1$ is ruled.  We conclude that
$M_1=\emptyset$ and the result follows from Lemma $1.1$ in \cite{DajGro}.\qed\vspace{1,5ex}

Observe that if  $f\colon\,M^n \rightarrow \R^N$  is a complete, simply connected parabolic submanifold, then $M^n$ is
diffeomorphic to $\R^n$ since its sectional curvature satisfies $K_M\leq 0$. In the ruled case, we have from Theorem
\ref{teo:hiper} that $M^n$ admits an isometric immersion as a ruled hypersurface with the same rulings. There are
many examples of complete  ruled hypersurfaces \cite{DajGro}. A simple example goes as follows: take  
$c\colon\,I\subset\R\rightarrow\R^{n+1}$ any unit speed curve, and let $E_0=dc/ds,E_1,\ldots, E_n$ a Frenet frame. It is easy to see that the hypersurface
$$
(s,t_1,\ldots,t_{n-1})\mapsto c(s)+\sum_{j=1}^{n-1} t_jE_{j+1}
$$
is complete.\vspace{1,5ex}

Given a nonruled  parabolic submanifold $f\colon\,M^n \rightarrow \R^N$ without Euclidean factor, let $\tilde{M}^n$
be the extension of  $f(M^n)$ in $\mathbb{R}^N$ obtained by extending each leaf of relative nullity of $f$ to a
complete affine Euclidean space in $\R^{n-2}$.
Our next and last result, describes the singular set of nonruled parabolic submanifolds without Euclidean factor
and dimension $n\geq 4$.

\begin{proposition}\po
Let $f\colon\,M^n \rightarrow \R^N$, $n\ge 4$, be a nonruled parabolic submanifold without Euclidean factor. Then the 
hypersurface given by
$$
\{\lambda \in \tilde{M}^n:\<\lambda,\xi_2^{s+1}\>=0\}
$$
is the singular set of $\tilde{M}$.
\end{proposition}

\proof Let $\Psi(\delta)=h(x)+ \delta$, $\delta \in \Lambda_s(x)$,
 be the parametrization in Theorem~\ref{theorem:polar}, where $h$ is any $s$--cross section of a polar surface $g$ to $f$. Without
loss of generality, we  assume that $h$ is a $\tau_*^g$-section. Being $(x,z)$  a coordinate system of $g$ with
$Z=\partial/\partial z$ asymptotic and $\{\eta_1,\ldots,\eta_k\}$ an orthonormal frame of $\Lambda_s$, we can
write
$$
\Psi(x,z,t_1,\ldots,t_k)=h(x,z)+ \sum_{j=1}^k t_j\eta_j(x,z)
$$
where $k=N-2s$ and $(t_1,\ldots,t_k)\in\R^k$. Recall that  $TM =\Lambda_{s-1}$ and $\Delta=\Lambda_s$. Thus, with
$X=\partial/\partial x$, we have that $\Psi(x,z,t_1,\ldots,t_k)$ is a singular point if and only if
$$
t_1(\nabla^\perp_X\eta_1)_{N_s}+t_2(\nabla^\perp_X\eta_2)_{N_s}\;\;\; \mbox{and}\;\;\;
t_1(\nabla^\perp_Z\eta_1)_{N_s}+t_2(\nabla^\perp_Z\eta_2)_{N_s}
$$
are linearly independents. By the definition of $\niu_s$, we have
$$
\<\nabla^\perp_Z\eta_1,\xi_2^s\>=\<\nabla^\perp_Z\eta_2,\xi_2^s\>=0.
$$
 Thus $t_1(\nabla^\perp_Z\eta_1)_{N_s}+t_2(\nabla^\perp_Z\eta_2)_{N_s}$ and
$\xi_2^s$ are normal fields. The above condition is now equivalent to $$
\< t_1(\nabla^\perp_X\eta_1)_{N_s} +
t_2(\nabla^\perp_X\eta_2)_{N_s},\xi_2^s\>=0 
$$ and, from Proposition \ref{Lemma1},  equivalent to
$$
\< t_1\eta_1+ t_2\eta_2,\xi_2^{s+1}\>=0.
$$
It follows that  $\lambda \in \tilde{M}^n$ is a singular point if and only if $ \<\lambda,\xi_2^{s+1}\>=0.$\qed

{\renewcommand{\baselinestretch}{1} \hspace*{-20ex}\begin{tabbing} \indent\= IMPA -- Estrada Dona
Castorina, 110
\indent\indent\=  Universidade da Beira Interior\\
\> 22460-320 -- Rio de Janeiro -- Brazil  \>
6201-001 -- Covilhã -- Portugal \\
\> E-mail: marcos@impa.br \> E-mail: pmorais@mat.ubi.pt
\end{tabbing}}

\end{document}